\documentclass[10pt]{article}
\usepackage{geometry}
\usepackage{hyperref,version}
\usepackage[T1]{fontenc}
\usepackage[utf8]{inputenc}
\usepackage[american]{babel}
\usepackage{amssymb}
\usepackage{amsmath,amsthm}
\numberwithin{equation}{section}
\usepackage{graphicx}
\usepackage{placeins}
\usepackage{xspace}
\usepackage{bm,bbm}
\usepackage{centernot}
\usepackage[numbers]{natbib}
\usepackage{enumerate,array}
\usepackage{xcolor}
\usepackage{dsfont}
\usepackage{enumitem}
\setlist{leftmargin=20pt}
\usepackage{Macros/mydef}
\usepackage{multirow}
\usepackage{diagbox}
\usepackage{lineno}
\usepackage{soul}
\usepackage{fullpage}

\usepackage{algorithm}
\usepackage{algpseudocode}
\usepackage{subcaption}
\usepackage{comment}
\usepackage{cancel}

\begin{document}
\newcommand\footnotemarkfromtitle[1]{%
\renewcommand{\thefootnote}{\fnsymbol{footnote}}%
\footnotemark[#1]%
\renewcommand{\thefootnote}{\arabic{footnote}}}

\title{Two nodally bound-preserving and mass-conservative discontinuous Galerkin methods for reaction-diffusion problems}

\author{Jie Shen\footnotemark[1]\, and Zuodong Wang\footnotemark[1]}%

\date{}

\maketitle

\renewcommand{\thefootnote}{\fnsymbol{footnote}}

\footnotetext[1]{School of Mathematical Sciences; Eastern Institute of Technology, Ningbo, Zhejiang, 315200, China. This work is supported in part by National Natural Science Foundation of China grants 12371409 and 12531015.}%

\renewcommand{\thefootnote}{\arabic{footnote}}
\begin{abstract}
We propose and analyze two high-order discontinuous Galerkin methods that preserve prescribed nodal bounds and total mass for reaction-diffusion problems with anisotropic and heterogeneous diffusion. Both methods combine a weighted interior-penalty prediction step with a convex correction in either the dG energy or the $L^2$ metric. We construct a pointwise bound-preserving and cellwise mass-preserving comparison function with optimal local approximation properties. This construction yields optimal dG-energy error estimates for the energy correction and, under a dual regularity assumption, optimal $L^2$ and dG-energy error estimates for the $L^2$ correction on quasi-uniform meshes. We further establish a general a posteriori estimate that combines a predictor estimator with the computable correction norm. The resulting estimators are globally reliable and efficient up to data oscillation, with parameter-dependent efficiency constants. Numerical experiments compare the two corrections and illustrate their use in adaptive computations.

\medskip
\noindent \textbf{Keywords.} 
Discontinuous Galerkin method, nodal bound-preservation, mass-conservation, prediction-correction method, reaction-diffusion problems, anisotropic and heterogeneous diffusion.

\medskip
\noindent \textbf{MSC.} 65N15, 65N30

\end{abstract}

\section{Introduction}

Numerical schemes with structure-preserving properties play an essential role in the numerical simulations in physics and engineering, especially in materials science, biology and fluid dynamics \cite{Gurtin_Polignone_Vinals_1996,Anderson_McFadden_Wheeler_1998,Fremond_2002,Vazquez_2007,Perthame_2015}. For example, positivity-preservation ensures that diffusion terms like $-\DIV(\sqrt{u}\nabla u)$ \cite{Vazquez_2007} are well-defined; mass-conservation is important for accurately reproducing long-time behavior of time-dependent problems \cite{Cheng_Shen_2022}; bound-preservation ensures that source terms such as $\log(1+u)-\log(1-u)$ \cite{Shen_Wang_2026} are well-defined. Within these models, many linear ones can be considered as linear reaction-diffusion problems, whereas nonlinear stationary and time-dependent problems can be formulated as linear reaction-diffusion problems after linearization or time-discretization. Another common feature of these models is the existence of sharp interfaces or singularities caused by heterogeneous and anisotropic diffusion tensor \cite{Ern_Stephansen_Zunino_2009}.

Based on the above discussion, we consider a class of linear reaction-diffusion PDEs (partial differential equations) which covers model features mentioned above: on a bounded, connected Lipschitz domain $\Dom\subset \Real^d$ with $d\in\{2,3\}$, let $u$ be the solution of
\begin{equation}\label{eq_model_pb}
    \mu u -\DIV (K\nabla u) = f \text{ in }\Dom,\qquad
    u = g_D \text{ on }\Gamma_D, \qquad K\nabla u \cdot \bn_{\Dom} = g_N \text{ on }\Gamma_N,
\end{equation}
where $K\in L^\infty(\Dom;\Real^{d\times d})$ is piecewise constant, symmetric, uniformly positive definite, anisotropic (the ratio between maximal and minimal eigenvalue of $K$ can be much larger than $1$) and heterogeneous ($K$ can have large jumps) diffusion matrix. $\mu>0$ is the reaction coefficient, $f\in L^2(\Dom)$, $|\Gamma_D\cap\Gamma_N|=0$, $\partial\Dom=\Gamma_D\cup\Gamma_N$, $g_D\in H^\frac{1}{2}(\Gamma_D)$ and $g_N\in L^2(\Gamma_N)$. 
With the help of trace operator $\gamma$, we introduce the trial and test function spaces 
\begin{equation}
    V\eqq \{ v\in H^1(\Dom)\ : \ \gamma(v)=g_D \text{ on } \Gamma_D \}, \qquad V_0\eqq \{ v\in H^1(\Dom)\ : \ \gamma(v)= 0  \text{ on } \Gamma_D \}.
\end{equation}
The model problem \eqref{eq_model_pb} can be rewritten in a weak formulation: find $u\in V$, such that
\begin{equation}\label{eq_model_pb_weak}
    a(u,v)=l(v) \quad \forall v\in V_0,
\end{equation}
where $a(v,w)=(\mu v,w)_\Dom +(K\nabla v,\nabla w)_\Dom$ for all $v,w\in H^1(\Dom)$ and $l(v)=(f,v)_\Dom +\int_{\Gamma_N}g_N v$ for all $v\in H^1(\Dom)$. 
We assume in addition that 
\begin{equation}\label{eq_maximum_principle}
    u(\xcoord)\in [a,b] \quad \text{ for almost all } \xcoord\in\Dom,
\end{equation}
with $-\infty<a<b<\infty$, and the total mass $m\eqq \int_\Dom u$ (hence $m/|\Dom|\in [a,b]$) is known a priori. We note that \eqref{eq_maximum_principle} is in general violated by standard numerical schemes, even with sufficiently fine mesh \cite{Herbin_Hubert_2008}.

The aim of this work is to find an efficient numerical method which is high-order, bound-preserving, mass-conservative, and able to capture and resolve the local features and singularities. For the finite element type schemes, it is practical to replace the pointwise bound-preservation by the nodal bound-preservation, hence in the rest of this work, the term \emph{bound-preservation} will hereafter refer to nodal bound-preservation.

Over the past decades, bound-preservation has remained a major challenge for finite element type schemes. In low-order settings (polynomial degree equals one), discrete maximum principles or bound-preserving properties can be achieved under suitable mesh conditions, such as the acute-angle condition for triangular meshes in two dimensions \cite{Barrenechea_John_Knobloch_2024}. In recent years, variational inequality formulations are applied to obtain high-order bound-preserving schemes in the pioneering works \cite{Barrenechea_Georgoulis_Pryer_Veeser_2024,Kirby_Shapero_2024}. Optimal convergence results in the energy norm have recently been established and extended to other model problems and spatial discretizations in \cite{Barrenechea_Georgoulis_Pryer_Veeser_2024,Amiri_Barrenechea_Pryer_2024,Amiri_Barrenechea_Pryer_2025,Dong_Ern_Wang_2025,Barrenechea_Pryer_Trenam_2025,Amiri_Barrenechea_Georgoulis_Pryer_2025,Shen_Wang_2026}. The optimal convergence rate in the $L^2$ norm has recently been established for piecewise affine finite elements in one dimension by \cite{Barrenechea_ChaumontFrelet_2026}. However, the above approaches generally do not ensure mass-conservation. This limitation has been addressed in a local sense in \cite{Barrenechea_Lederer_Rupp_2026} via an enriched Galerkin method with specifically designed iterative solver to address the resulting linear system with large condition number. To the best of our knowledge, in the literature, the following problems are important but remain open for high-order structure-preserving finite element type schemes: 1. How to efficiently ensure bound-preservation, mass-conservation and high-order accuracy simultaneously?  2. How to derive the error estimate in the $L^2$ norm in multidimensional setting? 3. How to derive a posteriori analysis and related adaptive scheme for singular solutions? 4. How to derive pointwise error estimates in the sense of \cite{Johnson_Schatz_Wahlbin_1987}? 

On the other hand, in the literature of time-dependent problems, a natural way to numerically ensure the physical constraints such as bound-preservation or mass-conservation is to
project (or correct) the numerical solutions onto the constrained admissible set, \eg the projection methods \cite{Bao_Du_2004,Hairer_Lubich_Wanner_2006,Bartels_Prohl_2007,Zhang_2017,Tong_Cai_2026}. We note that, for the bound-preserving and mass-conservative admissible set, the projection method is close to the Lagrange multiplier method \cite{Cheng_Shen_2022}. To the best of our knowledge, all methods mentioned above are mainly designed for problems with isotropic or mildly anisotropic diffusion, and focus on (quasi-) uniform mesh in space. Their technique will be adapted in this work for problems with anisotropic diffusion.

In order to capture the local feature and singular behavior, one natural idea is to consider the adaptive schemes with local refinement based on the a posteriori analysis (see, \eg \cite{Verfurth_2013} and the references therein). The challenging issues of minimal regularity and robustness with respect to $K$ of our model problem are investigated in \cite{Ern_Stephansen_Zunino_2009,Ern_Stephansen_Voralik_2010,Cai_He_Zhang_2017} and references therein. However, to the best of our knowledge, the bound-preserving high-order schemes are not analyzed in the a posteriori sense in the literature. Instead, the obstacle problems or the variational inequalities of the first kind are close to our model problem, and they are well studied in both a priori and a posteriori analysis (see \eg \cite{Veeser_2001,Burg_Schroder_2015,Han_2024} and references therein). These techniques can be adapted to our schemes, and the resulting estimators are compared with ours in \S\ref{sec_num_test}.

Overall, to the best of our knowledge, in the literature, there is no method which has simultaneously high-order accuracy, bound-preservation, mass-conservation, a priori and a posteriori error estimates, but these features are partially satisfied by different methods and analysis techniques from different communities. Based on the above literature review, we propose and analyze two prediction-correction type schemes to solve the model problem \eqref{eq_model_pb}, which satisfy all properties mentioned above, and one of them admits a priori $L^2$ error estimate. The prediction step \eqref{eq_scheme_step1} is given by standard discontinuous Galerkin (dG) approximation with harmonic weighting \cite{Dryja_2003}, and the correction steps \eqref{eq_scheme_step2} and \eqref{eq_scheme_step2_L2} are realized by projecting the predicted solution onto structure-preserving admissible set through a minimization problem in different metrics.

The main contributions of this work are summarized as follows:
\begin{itemize}
    \item \textbf{High-order structure-preserving schemes.} Two high-order bound-preserving and mass-conservative schemes are proposed and shown to be well-posed for the model problem \eqref{eq_model_pb}, based on the prediction-correction methodology. We obtain error estimates with similar parameter and mesh size dependencies as in the literature, meanwhile additional structure-preserving properties are satisfied.

    \item \textbf{Structure-preserving dG candidate with optimal approximation properties.} Based on a cellwise convex combination of $L^2$ projection and its mean value, we construct a (pointwise) bound-preserving and cellwise mass-conservative dG function $\Pi_h u$ which converges to $u\in \calC(\bar{\Dom};[a,b])\cap H^1(\Dom)$ as $h\to 0$, and optimal convergence rate is established in $H^k$ norm, $k=0,...,p+1$ and $L^\infty$ norm, for $u$ in $H^{p+1}(\calT_h)$ and $W^{p+1,\infty}(\calT_h)$, respectively. This may be of independent interest, since it provides an optimal candidate for error analysis of other structure-preserving schemes. 
    
    \item \textbf{Scheme with dG energy correction.} For the scheme \eqref{eq_scheme_step1}-\eqref{eq_scheme_step2}, we prove, for the dG energy norm, the optimal a priori  convergence rate, under standard regularity assumptions in Theorem \ref{thm_apriori_energy}, without quasi-uniform mesh and $W^{p,\infty}$-regularity assumptions used in our previous work \cite{Shen_Wang_2026}; for the a posteriori analysis, we use the triangle inequality and the a posteriori approximation properties of the prediction step to construct error estimators \eqref{eq_eta_energy} and establish the global reliability and efficiency of these estimators in Theorem \ref{thm_aposteriori}.
    
    \item \textbf{Scheme with $L^2$ correction.} For the scheme  \eqref{eq_scheme_step1}-\eqref{eq_scheme_step2_L2}, we prove the optimal a priori convergence rate for the $L^2$ norm under standard regularity assumption, and the optimal a priori convergence rate for dG energy norm requires in addition the quasi-uniformity on meshes, in Theorem \ref{Thm_apriori_L2}. A posteriori analysis is also derived in Theorem \ref{thm_aposteriori_L2}, with the same technique used in the dG correction scheme. We also note that our a posteriori analysis framework can be extended to other correction steps. See Corollary \ref{cor_arbitrary_correction} and Remark \ref{rmk_other_correction} for details.

    \item \textbf{Numerical comparison.} Owing to the block-diagonal structure of mass matrix of dG method, the scheme with $L^2$ correction is computationally less expensive than scheme with dG energy correction. However, we observe in \S\ref{sec_num_test} that, the correction realized in the $L^2$ metric is less smooth than the correction realized by the dG energy, hence the usage of two schemes depends on the resolution of the diffusion. In addition, the comparison between our estimators and estimators derived based on variational inequality techniques are also conducted.
    
\end{itemize}

The paper is organized as follows. \S\ref{sec_scheme} presents the discrete setting, the prediction-correction schemes and their fundamental properties. \S\ref{sec_apriori} develops a priori error analysis of the scheme with dG energy correction. In \S\ref{sec_aposteriori}, we establish the a posteriori error analysis for the scheme with dG energy correction. In \S\ref{sec_L2_correction}, we extend and adapt the argument in \S\ref{sec_apriori} and \S\ref{sec_aposteriori} to the scheme with $L^2$ correction.
Numerical validation and comparison of two schemes are given in \S\ref{sec_num_test}, followed by concluding remarks in \S\ref{sec_conclusion}.

Throughout the paper, we use standard notation for Sobolev spaces. The $L^2(S)$ inner product is written as $(\cdot,\cdot)_S$ with associated norm $\norm{\cdot}_S$ for any $S\subset \Dom$. Generic constants $C>0$ may change from line to line but are always independent of the mesh size, $\mu$, $K$, $f$, $u$, $g_D$ and $g_N$, but may depend on mesh shape regularity, domain $\Dom$, polynomial degree and dG stabilization parameter. For brevity, spatial dependence is often suppressed when clear from context (\eg $K=K(\xcoord)$).

\section{Discrete setting and numerical scheme}\label{sec_scheme}
In this section, we present the discrete setting, two prediction-correction high-order structure-preserving schemes, their fundamental properties and a practical solver for the correction step.

\subsection{Discrete setting}

For the spatial discretization, we consider a simplicial conforming mesh $\{\calT_h\}_{h>0}$ of $\Dom$.  A generic mesh cell is denoted by $T\in \calT_h$. For simplicity, we assume that $K$ is constant on each cell $T\in\calT_h$. We define the mesh size parameter $h\eqq \max_{T\in\calT_h}h_T $ and $ h_T\eqq \text{diam}(T)$. We also define $h_{\min} \eqq \min_{T\in\calT_h} h_T$. The set of mesh faces is defined by $\calF\eqq \calF_D\cup\calF_N\cup\calF_i$, with $\calF_D$, $\calF_N$ and $\calF_i$ representing the sets of Dirichlet, Neumann and interior faces, respectively. For each $F\in\calF$, we define $h_F\eqq \text{diam}(F)$ as the face diameter. We also set $\calT_F$ as the set of cells $T$ such that $|\partial T \cap F|\neq 0$, and $\Tilde{\calF}_T$ as the set of faces sharing a vertex with $T$. For each face $F\in\calF_i$, we define $\bn_F$ as the globally fixed unit normal vector on $F$. For $F\in\calF_i$, it points from $T_+$ to $T_-$, for some $T_\pm\in\calT_h$ with $|\partial T_\pm \cap F| >0$, and for $F\in\calF_D\cup\calF_N$, it is the outward unit normal vector of the domain. We assume only that the mesh is shape-regular, rather than the \emph{acute angle} assumption which guarantees the discrete maximum principle in the literature (see for example \cite{Barrenechea_John_Knobloch_2024}). On each cell $T\in\calT_h$, we set $K_{\max,T}$ and $K_{\min,T}$ as the maximum and minimum eigenvalues of $K$ on $T\in\calT_h$, respectively. We also define $K_{\max}\eqq \max_{T\in\calT_h} K_{\max,T}$ and $K_{\min}\eqq \min_{T\in\calT_h} K_{\min,T}$.

We consider the dG space with a Lagrange basis of polynomial degree $p\geq 1$. The global shape functions are
denoted by $\{\varphi_i\}_{i \in \DoFSet}$. The corresponding dG space is denoted by 
\begin{equation}
    V_h^p\eqq\text{span}\{\varphi_i\}_{i \in \DoFSet}= \left\{ v\in L^\infty(\Omega)\ : \  v|_T \in \mathbb{P}_p(T)\ \forall T\in\calT_h \right\},
\end{equation}
where $\mathbb{P}_p(T)$ denotes the local Lagrange finite element space of degree $p$ over $T$. The set of physical nodes of associated Lagrange basis (nodes with same physical coordinate are considered as the same node) is denoted by $\calV^{cG}$.

We denote the broken Sobolev spaces by 
\begin{equation}
    W^{k,q}(\calT_h) \eqq \left\{  v\in L^q(\Dom) \ : \ v|_T \in W^{k,q}(T)\ \forall T\in\calT_h  \right\},\quad H^{s}(\calT_h) \eqq W^{s,2}(\calT_h).
\end{equation}
$\nabla_h$, $\nabla_h\cdot $, $\nabla_h^2$ and $\Delta_h$ are the cellwise gradient, divergence, Hessian, and Laplacian, respectively, \ie $\nabla_h|_T=\nabla$, $\nabla_h\cdot|_T=\DIV$, $\nabla_h^2|_T=\nabla^2$ and $\Delta_h|_T=\Delta$ for all $T\in\calT_h$.

\subsection{Numerical schemes with dG energy correction}
We consider the dG method with harmonic weighting \cite{Dryja_2003,Ern_Stephansen_Zunino_2009,DiPietro_Ern_2012} to approximate the model problem \eqref{eq_model_pb}, as it is more robust with respect to jumps of $K$ than the standard continuous finite element method or the unweighted symmetric interior penalty dG.
More precisely, the prediction step is given by finding $\Bar{u}_h\in V_h^p$ such that
\begin{equation}\label{eq_scheme_step1}
    a_h(\Bar{u}_h,v_h) = l_h(v_h) \quad \forall v_h\in V_h^p,
\end{equation}
where
\begin{subequations}
\begin{align}
    a_h(v,w) &\eqq \int_\Dom \left( \mu vw+ K \nabla_h v \cdot \nabla_h w \right) \\
    &\quad + \sum_{F\in\calF_i\cup \calF_D} \int_F \left( -\avg{K\nabla_h v}\cdot \bn_F \jump{w} - \avg{K\nabla_h w}\cdot \bn_F \jump{v} + \gamma_0 \frac{\sigma_F}{h_F} \jump{v}\jump{w} \right),\nonumber\\
    l_h(v) &\eqq (f,v) + \int_{\Gamma_N}g_N v + \sum_{F\in\calF_D} \int_F \left( -K\nabla_h v \cdot \bn_F g_D + \gamma_0 \frac{\sigma_F}{h_F} v g_D\right),\\
    \avg{v} |_F &\eqq \omega_F^+ v^+ + \omega_F^-v^- \quad \forall F\in\calF_i, \qquad \avg{v} |_F = v^+ \quad \forall F\in\calF_D\cup\calF_N,\\
    \jump{v} |_F &\eqq v^+ - v^- \quad \forall F\in\calF_i, \qquad \jump{v} |_F = v^+ \quad \forall F\in\calF_D\cup\calF_N,
\end{align}
\end{subequations}
and $\omega_F^\pm\eqq \frac{\delta_F^\mp}{\delta_F^++\delta_F^-}$, $\sigma_F \eqq \frac{\delta_F^+\delta_F^-}{\delta_F^++\delta_F^-}$, $\delta_F^\pm\eqq \bn_F^T K^\pm \bn_F$, $v^\pm = v|_{T_\pm}$ for $F\in\calF_i$, $\sigma_F\eqq \delta_F \eqq \bn_F^T (K|_{T_+}) \bn_F$, $v^+=v|_{T_+}$ for $F\in\calF_D\cup\calF_N$ with $T_+$ the cell satisfying $F\subset \partial T_+$. $\gamma_0 > 0$ is a user-defined stabilization parameter. We observe that, for the weighted averaging operator, we have, for all $F\in\calF_i$, and all $v,w\in V_h^p + H^1(\Dom)$,
\begin{equation}
    \jump{vw}=\jump{v}\avg{w} + \jump{w}\avgbar{v}=\jump{v}\avgbar{w} + \jump{w}\avg{v},
\end{equation}
with $\avgbar{v}=\omega_F^-v^+ + \omega_F^+ v^-$.

Then, the structure-preserving numerical solution is constructed by the correction step 
\begin{equation}\label{eq_scheme_step2}
    u_h=\text{argmin}_{v_h\in X_h} \frac{1}{2}a_h(v_h-\Bar{u}_h,v_h-\Bar{u}_h),
\end{equation}
where $X_h\subset V_h^p$ is an abstract nonempty, convex and compact structure-preserving subset of $V_h^p$. In this work, we consider only a bound-preserving subset $X_h^{a,b}$ defined as
\begin{equation}\label{eq_Xh_ab}
    X_h^{a,b}\eqq \left\{ v_h=\sum_{i\in\DoFSet} V_i\varphi_i \ : \ V_i\in[a,b] \quad \forall i\in\DoFSet \right\},
\end{equation}
and a bound-preserving and mass-conservative subset $X_h^{a,b,m}$ given by
\begin{equation}\label{eq_Xh_abm}
    X_h^{a,b,m}\eqq \left\{ v_h=\sum_{i\in\DoFSet} V_i\varphi_i \ : \ V_i\in[a,b] \quad \forall i\in\DoFSet, \int_\Dom v_h=m \right\}.
\end{equation}
Other choices of structure-preserving subsets of $V_h^p$ are also possible, although their practical usefulness depends on whether the resulting minimization problem can be solved efficiently.

\subsection{Numerical scheme with $L^2$ correction}
In this part, we replace the correction step \eqref{eq_scheme_step2} by its $L^2$ counterpart
\begin{equation}\label{eq_scheme_step2_L2}
    \mathfrak{u}_h=\text{argmin}_{v_h\in X_h} \frac{1}{2} \norm{v_h-\bar{u}_h}_\Dom^2.
\end{equation}
This variant is easier to solve using the PDAS solver, compared with the original correction step \eqref{eq_scheme_step2}. Indeed, since the mass matrix of the current $L^2$ functional is block diagonal, its corresponding linear system can be solved blockwise, unlike other variational inequality based structure-preserving schemes \cite{Barrenechea_Georgoulis_Pryer_Veeser_2024,Kirby_Shapero_2024,Dong_Ern_Wang_2025,Shen_Wang_2026,Barrenechea_Lederer_Rupp_2026}. Therefore, the computational cost of PDAS solver based on this correction step is negligible compared with the prediction step \eqref{eq_scheme_step1} or the dG energy correction step \eqref{eq_scheme_step2}. In addition, this variant leads to optimal error estimate in the $L^2$ norm, owing to the trivial triangle inequality and approximation properties of the minimization problem \eqref{eq_scheme_step2_L2}:
\begin{equation}
    \norm{u-\mathfrak{u}_h}_{\Dom} \leq \norm{u-\bar{u}_h}_{\Dom} + \norm{\bar{u}_h-\mathfrak{u}_h}_{\Dom} \leq 2 \norm{u-\bar{u}_h}_{\Dom} + \inf_{v_h\in X_h} \norm{v_h-u}_{\Dom}.
\end{equation}
Moreover, the optimal error estimate in the dG energy norm is provable for quasi-uniform meshes, with the help of discrete trace and inverse inequalities \eqref{eq_dis_trace_inv}. These a priori estimates will be established in \S\ref{sec_L2_correction}.

To the best of our knowledge, this is the first time that the optimality of error estimates in both $L^2$  and energy norms has been proved for high-order structure-preserving methods in multi-dimensional setting. This cheap correction step is, however, less accurate for underresolved anisotropic diffusion, compared with the dG energy based correction step \eqref{eq_scheme_step2}, we refer to \S\ref{sec_num_test} for numerical evidence.

Moreover, the scheme with $L^2$ correction \eqref{eq_scheme_step1}-\eqref{eq_scheme_step2_L2} can be considered as a variant of some existing projection methods, based on the following remark:
\begin{remark}[Relations to existing methods]
    The Lagrange multiplier method with Fourier spectral discretization \cite{Cheng_Shen_2022} and the projection method with finite difference discretization \cite{Tong_Cai_2026} are in fact close to our scheme with $L^2$ correction \eqref{eq_scheme_step1}-\eqref{eq_scheme_step2_L2}. Indeed, for Fourier spectral method or finite difference method with structured mesh (grid), their correction step is given by minimizing the functional $h^d\sum_{i\in\DoFSet} (V_i-\bar{U}_i)^2$, which coincides with the projection induced by the corresponding discrete $L^2$ inner product.
\end{remark}

\subsection{Fundamental properties}
We prove in this section the well-posedness of schemes \eqref{eq_scheme_step1}-\eqref{eq_scheme_step2} and \eqref{eq_scheme_step1}-\eqref{eq_scheme_step2_L2} and show their relation with variational inequalities.

We recall the discrete trace and inverse inequalities which will be frequently used in this work:
\begin{lemma}[Discrete trace and inverse inequalities]
    Let $v_h\in V_h^p$. For any $F\in \partial T$ of some $T\in\calT_h$, we have
    \begin{equation}\label{eq_dis_trace_inv}
        \norm{v_h}_F \leq C h_F^{-1/2} \norm{v_h}_T ,\quad \norm{\nabla v_h}_T \leq C h_T^{-1}\norm{v_h}_T.
    \end{equation}
\end{lemma}

\begin{proof}
    The proof can be found in \cite[Section 1.4.3]{DiPietro_Ern_2012}.
\end{proof}

Continuous trace inequality is also needed:
\begin{lemma}[Continuous trace inequality]
    Let $v\in H^1(T)$ for $T\in\calT_h$. We have
    \begin{equation}\label{eq_con_trace}
        \norm{v}_{\partial T}^2 \leq C (\norm{\nabla v}_T + h_T^{-1} \norm{v}_T)\norm{v}_T.
    \end{equation}
\end{lemma}

\begin{proof}
    The proof can be found in \cite[Lemma 1.49]{DiPietro_Ern_2012}.
\end{proof}

We introduce the dG norm as follows:
\begin{equation}\label{eq_dG_norm}
    \dgnorm{v}^2 \eqq \sum_{T\in\calT_h} \norm{\mu^{\frac{1}{2}} v }_T^2 + \sum_{T\in\calT_h} \norm{K^\frac{1}{2}\nabla v}_T^2 + \sum_{F\in \calF_D\cup \calF_i } \frac{\sigma_F}{h_F}\norm{\jump{v}}_F^2\qquad \forall v\in V_h^p+H^1(\Dom).
\end{equation}

Using the dG norm, we have the well-posedness of the prediction step \eqref{eq_scheme_step1}:
\begin{lemma}[Well-posedness of prediction step]
    For sufficiently large $\gamma_0$, the bilinear form $a_h$ is coercive and bounded in the sense
    \begin{equation}\label{eq_coercivity}
        a_h(v_h,v_h) \geq C \dgnorm{v_h}^2, \quad |a_h(v_h,w_h)|\leq C\dgnorm{v_h}\dgnorm{w_h} \qquad \forall v_h,w_h\in V_h^p,
    \end{equation}
    and the prediction step \eqref{eq_scheme_step1} is well-posed.
\end{lemma}
\begin{proof}
    The proof is based on the standard Lax-Milgram theorem. We refer to \cite[Lemma 3.1]{Ern_Stephansen_Zunino_2009} for details.
    
\end{proof}

\begin{lemma}[Well-posedness of dG energy correction step]
    For sufficiently large $\gamma_0$, the correction step \eqref{eq_scheme_step2} is well-posed.
\end{lemma}
\begin{proof}
    Since $a_h$ is symmetric and coercive, the quadratic functional in \eqref{eq_scheme_step2} is strictly convex, hence the optimization problem \eqref{eq_scheme_step2} has a unique solution if its admissible set $X_h$ is nonempty, compact and convex. For $X_h=X_h^{a,b,m}$, since $m/|\Dom|\in [a,b]$ is an admissible element in $X_h$, the admissible set is nonempty. The compactness of $X_h$ follows from the fact that $V_h^p$ has finite dimension and $X_h$ is closed and bounded in $V_h^p$; for any $v_h\eqq\sum_{i\in\DoFSet} V_i\varphi_i  ,w_h\eqq\sum_{i\in\DoFSet} W_i\varphi_i \in X_h$ and $\theta\in[0,1]$, we have $\theta V_i + (1-\theta) W_i \in [a,b]$ and $\int_\Dom \theta v_h + (1-\theta)w_h=m$ for $\int_\Dom v_h =\int_\Dom w_h=m$, which leads to the convexity of $X_h^{a,b,m}$. The same argument also leads to the nonemptiness, convexity and compactness of the case $X_h=X_h^{a,b}$.
\end{proof}

Following the same argument, we immediately obtain the well-posedness of the $L^2$ correction step:
\begin{lemma}[Well-posedness of $L^2$ correction step]
    The correction step \eqref{eq_scheme_step2_L2} is well-posed.
\end{lemma}

Next, we prove that the proposed prediction-correction step is actually equivalent to a variational inequality:
\begin{lemma}[Variational inequality representation of dG correction based scheme]
    The solution to \eqref{eq_scheme_step1}-\eqref{eq_scheme_step2} is equivalently characterized by the following variational inequality:
    \begin{equation}\label{eq_VI}
        a_h(u_h,v_h-u_h) \geq l_h(v_h-u_h) \qquad \forall v_h\in X_h.
    \end{equation}
    In addition, it satisfies the equation
    \begin{equation}\label{eq_LM}
        a_h(u_h,v_h) = l_h(v_h) + (\sigma_h,v_h) \eqq \Tilde{l}_h(v_h) \qquad \forall v_h\in V_h^p,
    \end{equation}
    with
    \begin{equation}\label{eq_sigma_h}
        \sigma_h \eqq \begin{cases}
            \lambda_a-\lambda_b \text{ if } X_h=X_h^{a,b}\\
            \lambda_a-\lambda_b-\xi_h \text{ if } X_h=X_h^{a,b,m}
        \end{cases},
    \end{equation}
    where $\lambda_a$, $\lambda_b$, $\xi_h$ are the dG function representation of Lagrange multipliers for lower bound, upper bound and mass conservation, respectively.
\end{lemma}

\begin{remark}[Variational inequality]
    Since such prediction-correction scheme is equivalent to a variational inequality, our scheme can be considered as an alternative solver for variational inequalities like obstacle problems \cite{Han_2024}, and our analysis in the next sections can be applied directly.
\end{remark}

\begin{proof}
    We prove \eqref{eq_VI} first, then derive \eqref{eq_LM} and \eqref{eq_sigma_h} by invoking the KKT condition.

    For all $v_h\in X_h$, we expand $a_h(v_h-\Bar{u}_h,v_h-\Bar{u}_h)$ to find
    \begin{align*}
        a_h(v_h-\Bar{u}_h,v_h-\Bar{u}_h) = a_h(v_h,v_h) -2 a_h(v_h,\Bar{u}_h) + a_h(\Bar{u}_h,\Bar{u}_h).
    \end{align*}
    Hence \eqref{eq_scheme_step2} is equivalent to
    \begin{equation*}
        u_h=\text{argmin}_{v_h\in X_h} \frac{1}{2}a_h(v_h,v_h) - a_h(v_h,\Bar{u}_h).
    \end{equation*}
    Since this optimization problem has a unique minimizer, it is equivalent to its first-order condition
    \begin{equation*}
        a_h(u_h,v_h-u_h) \geq a_h(\Bar{u}_h,v_h-u_h) = l_h(v_h-u_h) \qquad v_h \in X_h,
    \end{equation*}
    with the help of \eqref{eq_scheme_step1}. 

    To derive the variational equation \eqref{eq_LM}, we only present the proof of the case $X_h=X_h^{a,b,m}$ here, since the other case can be deduced by the same argument. We notice that the KKT conditions of the equivalent constrained minimization problem \eqref{eq_scheme_step2} are
    \begin{align*}
        &a_h(u_h,v_h) = l_h(v_h) - \xi_h \int_\Dom v_h - \langle \mu_b-\mu_a , v_h \rangle  \qquad \forall v_h\in V_h^p,\\
        &\int_\Dom u_h = m, \qquad \mu_{b,i}(U_i-b) = 0, \quad \mu_{a,i}(a-U_i) = 0, \quad \mu_{b,i},\mu_{a,i} \geq 0\quad \forall i\in \DoFSet,
    \end{align*}
    where $\langle \cdot,\cdot\rangle$ is the dof-wise Euclidean inner product, $\xi_h\in\Real$, $\mu_b\eqq (\mu_{b,1},...,\mu_{b,|\DoFSet|})$ and $\mu_a\eqq (\mu_{a,1},...,\mu_{a,|\DoFSet|})$ are Lagrange multipliers for mass-conservation, upper and lower bound constraints, respectively. In order to define dG functions $\lambda_b,\lambda_a\in V_h^p$, we introduce the Riesz representation
    \begin{equation*}
        (\lambda_c,v_h)_\Dom=\sum_{i\in\DoFSet} \mu_{c,i}V_i, \qquad \forall v_h\eqq \sum_{i\in\DoFSet} V_i\varphi_i, \quad c\in\{a,b\}.
    \end{equation*}
    Finally, we set $\sigma_h\eqq \lambda_a - \lambda_b -\xi_h$ to conclude the proof.
\end{proof}

Then, for the numerical scheme with $L^2$ correction \eqref{eq_scheme_step1}-\eqref{eq_scheme_step2_L2}, we can introduce a metric lifting $\calL_h:V_h^p\to V_h^p$ defined by
\begin{equation}
    (\calL_h v_h,w_h)_{\Dom} = a_h(v_h,w_h) \quad \forall v_h,w_h\in V_h^p.
\end{equation}
Since $\mathfrak{u}_h$ is the minimizer of \eqref{eq_scheme_step2_L2}, there exists a Lagrange multiplier $\sigma_{h,0}\in V_h^p$ such that
\begin{equation}
    \mathfrak{u}_h-\bar{u}_h = \sigma_{h,0}.
\end{equation}
Then we apply the metric lifting operator on it to find
\begin{equation}\label{eq_scheme_L2_LM}
    a_h(\mathfrak{u}_h,v_h) =l_h(v_h) + (\calL_h\sigma_{h,0},v_h) \quad \forall v_h \in V_h^p,
\end{equation}
which has a similar form to \eqref{eq_LM}. This reformulation is useful in the a posteriori analysis of variational inequalities \cite{Veeser_2001,Braess_2005}. We postpone the discussion of this point later in \S\ref{sec_num_test}.

\subsection{Practical implementation}
Although the minimization problem \eqref{eq_scheme_step2} or \eqref{eq_scheme_step2_L2} is a standard quadratic program, standard optimization solvers such as OSQP \cite{Stellato_Banjac_Goulart_Bemporad_Boyd_2020} are numerically observed to be expensive. Therefore, we present a solver based on the primal-dual active-set (PDAS) method \cite{Hintermuller_Ito_Kunisch_2002} here, and also give an explicit construction of $\sigma_h$. We also note that, PDAS is not the only efficient solver, Lagrange multiplier based solvers \cite{Tong_Cai_2026,Cheng_Shen_2022} are also widely used.

The basic idea of the PDAS method is to fix the bound-violating degrees of freedom at the corresponding bounds and determine the remaining degrees of freedom by solving the resulting linear system. This process is repeated until the bound-violating dof set no longer changes. The Lagrange multipliers are reconstructed from the KKT conditions of the variational inequality \eqref{eq_VI}. 
Such PDAS method for $X_h=X_h^{a,b,m}$ is given in Algorithm \ref{alg_BPMC}.
PDAS solver for $X_h^{a,b}$ is well-established in \cite{Hintermuller_Ito_Kunisch_2002}, hence we omit it for brevity. 

\begin{algorithm}[!htbp]
\caption{Bound-preserving and mass-conservative PDAS}\label{alg_BPMC}
\begin{algorithmic}

\Require System matrix $A\in\Real^{|\DoFSet|\times |\DoFSet|}$, prediction solution $\Bar{U}\in\Real^{|\DoFSet|}$, lumped mass vector $M_L=(\int_\Dom \varphi_1,...,\int_\Dom \varphi_{|\DoFSet|})\eqq (m_1,...,m_{|\DoFSet|})\in\Real^{|\DoFSet|}$, mass matrix $M\in\Real^{|\DoFSet|\times |\DoFSet|}$, maximal iteration number $k_{\max}$.
\Ensure Corrected solution $U$, Lagrange multiplier $\sigma_h$.

\State Set the target mass $M_0 \eqq m$, local mass $m_S=(m_i)_{i\in S}$ for any $S\subset \{1,...,|\DoFSet|\}$, initialize $U^0\eqq \bar U$, $\ell^0\eqq 0$, and choose
\[
\mathcal A_a^0 \eqq  \{i:\bar U_i\le a\},\qquad
\mathcal A_b^0 \eqq  \{i:\bar U_i\ge b\},\qquad \mathcal F^0 \eqq  \{1,\ldots,|\DoFSet|\}\setminus
(\mathcal A_a^0\cup\mathcal A_b^0)
\]

\For{$k=0,1,2,\ldots, k_{\max}$}
\State Define the active part $U_A^k$ by $(U_A^k)_i =
\begin{cases}
a, & i\in \mathcal A_a^k,\\
b, & i\in \mathcal A_b^k,\\
0, & i\in \mathcal F^k
\end{cases}$

\State Compute $M_A^k \eqq  \sum_{i\in\mathcal A_a^k}m_i a + \sum_{i\in\mathcal A_b^k}m_i b, \qquad d_F^k \eqq  M_0-M_A^k$

\If{$m_{\mathcal F^k}=0$}
    \If{$d_F^k=0$}
        \State Set $U^{k+1}\eqq U_A^k$ and
        $\ell^{k+1}\eqq\ell^k$

        \If{$\mathcal F^k\neq\emptyset$}
            \State Let
            $\mathcal A^k\eqq
            \mathcal A_a^k\cup\mathcal A_b^k$
            \State Solve $A_{\mathcal F^k\mathcal F^k}
                U_{\mathcal F^k}^{k+1}
                =
                (A\bar U)_{\mathcal F^k}
                -
                A_{\mathcal F^k\mathcal A^k}
                (U_A^k)_{\mathcal A^k}$
        \EndIf
    \Else
        \State Choose
        $j\in\mathcal A_a^k$ with $d_F^k m_j>0$,
        or $j\in\mathcal A_b^k$ with $d_F^k m_j<0$

        \State Set
        $\mathcal A_a^k\gets\mathcal A_a^k\setminus\{j\}$,
        $\mathcal A_b^k\gets\mathcal A_b^k\setminus\{j\}$,
        and
        $\mathcal F^k\gets\mathcal F^k\cup\{j\}$

        \State \textbf{restart} the current iteration
        from the construction of $U_A^k$
    \EndIf
\Else
    \State Let $\mathcal A^k\eqq \mathcal A_a^k\cup\mathcal A_b^k$, solve $A_{\mathcal F^k\mathcal F^k} Z_{\mathcal F^k}
    =
    (A\bar U)_{\mathcal F^k}
    -
    A_{\mathcal F^k\mathcal A^k}(U_A^k)_{\mathcal A^k}$ and $A_{\mathcal F^k\mathcal F^k} W_{\mathcal F^k}
    =
    m_{\mathcal F^k}$

    \State Compute $\ell^{k+1}
    =
    \frac{
    m_{\mathcal F^k}^T Z_{\mathcal F^k}-d_F^k
    }{
    m_{\mathcal F^k}^T W_{\mathcal F^k}
    }$
    \State update $ U_{\mathcal F^k}^{k+1}
    =
    Z_{\mathcal F^k}-\ell^{k+1}W_{\mathcal F^k}$, $U_i^{k+1}=a,\ i\in\mathcal A_a^k$, $U_i^{k+1}=b,\ i\in\mathcal A_b^k$
\EndIf

\State Compute the KKT gradient $Y^{k+1}\eqq A(U^{k+1}-\bar U)+\ell^{k+1}M_L$

\State Update 
$
\mathcal A_a^{k+1}
\eqq 
\left\{
i:\ U_i^{k+1}-\frac{1}{|A_{ii}|}Y_i^{k+1}\le a
\right\}$,
$
\mathcal A_b^{k+1}
\eqq 
\left\{
i:\ U_i^{k+1}-\frac{1}{|A_{ii}|}Y_i^{k+1}\ge b
\right\} $,
and
$
\mathcal F^{k+1}
\eqq 
\{1,\ldots,|\DoFSet|\}\setminus
(\mathcal A_a^{k+1}\cup\mathcal A_b^{k+1})$

\If{$\mathcal A_a^{k+1}=\mathcal A_a^k$ and $\mathcal A_b^{k+1}=\mathcal A_b^k$}
    \State Set $U\eqq U^{k+1}$ and $\ell\eqq \ell^{k+1}$
    \State \textbf{break}
\EndIf

\EndFor

\State Recompute $\sigma_{h,i}\eqq (M^{-1}A(U-\bar U))_i$ for all $i\in\DoFSet$

\State \textbf{return} $U$, $\sigma_h$

\end{algorithmic}
\end{algorithm}

\section{A priori analysis}\label{sec_apriori}

We assume in this section that $u\in H^2(\calT_h)$ to simplify the presentation and to avoid additional techniques required for less regular solutions, e.g., \cite{Ern_Guermond_2022}. Moreover, our argument can be smoothly extended to less regular case, provided that convergence of the prediction step is ensured.

In the first part, we introduce a pointwise bound-preserving and cellwise mass-conservative dG candidate with optimal convergence rate in standard Sobolev norms; in the second part, we derive the Céa-type estimate over structure-preserving set $X_h$, and use the optimal candidate constructed in the first part to derive the optimal convergence of the numerical solution delivered by the prediction-correction scheme \eqref{eq_scheme_step1}-\eqref{eq_scheme_step2}.

\subsection{Structure-preserving candidate with optimal approximation property}
The first step of a priori error estimates is to construct a structure-preserving candidate which has optimal convergence properties. Indeed, the Céa-type estimate is formulated over $X_h$, so the comparison function must belong to this set. Standard elliptic projections and Lagrange interpolants need not satisfy the constraints; we therefore construct a candidate in $X_h$ with optimal approximation properties.

Before constructing such a function, we first present a technical lemma for mean-free polynomials:
\begin{lemma}[$L^\infty$ norm for mean-free polynomials]\label{lemma_poly_infty_min_max}
    For any $T\in\calT_h$, let $v_h$ be a non-constant polynomial on $T$ of order $p$, with $\int_T v_h=0$. We have
    \begin{equation}
        \norm{v_h}_{L^\infty(T)} \leq C \max_{\xcoord\in \bar{T}}v_h(\xcoord), \qquad \norm{v_h}_{L^\infty(T)} \leq C \Big( -\min_{\xcoord\in \bar{T} }v_h(\xcoord) \Big).
    \end{equation}
\end{lemma}

\begin{proof}
    Based on the norm equivalence in finite dimensional space, it is proved in \cite[Appendix C]{Zhang_2017} that
    \begin{equation*}
        \frac{\norm{v_h}_{L^\infty(T)}}{\max_{\xcoord\in \bar{T}}v_h(\xcoord)} \leq C,
    \end{equation*}
    which is the first inequality we want. Then, we replace $v_h$ by $-v_h$ in the above inequality to conclude the result.
\end{proof}

Using the preceding lemma, we can construct an optimal structure-preserving candidate, in the following lemma:

\begin{lemma}[Optimal candidate in $X_h^{a,b,m}$]\label{lemma_opt_vh_abm}
    Let $v\in \calC(\bar{\Dom};[a,b])\cap H^1(\Dom)$ with total mass $\int_\Dom v = m$. There exists a function $\Pi_h v\in X_h^{a,b,m}$, such that it satisfies
    \begin{equation}
        a\leq (\Pi_h v)|_T \leq b , \qquad \int_T \Pi_h v =\int_T v \quad \forall T\in\calT_h,
    \end{equation}
    and the convergence property
    \begin{equation}
        \norm{v-\Pi_h v}_{L^\infty(\Dom)},\norm{v-\Pi_h v}_{L^2(\Dom)},\norm{\nabla_h(v-\Pi_h v)}_{L^2(\Dom)} \to 0 \text{ as } h\to 0.
    \end{equation}
    Moreover, in each $T\in\calT_h$, we have
    \begin{align}
        \norm{\nabla^k(v-\Pi_h v)}_T &\leq C  h_T^{p+1-k} \norm{\nabla^{p+1} v}_T \quad\text{ if } v|_T\in H^{p+1}(T) \quad \forall k \in \{0,1,...,p+1\},\\
        \norm{v-\Pi_h v}_{L^\infty(T)} &\leq C h_T^{p+1} \norm{\nabla^{p+1} v}_{L^\infty(T)} \quad\text{ if }v|_T\in W^{p+1,\infty}(T).
    \end{align}
\end{lemma}

\begin{proof}
    Let $\Tilde{v}_h \eqq \Pi_h^{L^2} v$ be the $L^2$ projection. Such projection has stability and optimal approximation properties in $L^2$, $L^\infty$ , and broken $H^k$-norms for $k\in\{1,...,p+1\}$, and it is locally mass-conservative \cite[Section 11.5.3]{Ern_Guermond_FEM1_2021}. On each cell $T\in\calT_h$, we define $\bar{v}_T\eqq \frac{1}{|T|}\int_T \Tilde{v}_h=\frac{1}{|T|}\int_T v\in[a,b]$, $M_T \eqq \max_{\xcoord\in \bar{T}} \Tilde{v}_h$ and $m_T \eqq \min_{\xcoord\in \bar{T}} \Tilde{v}_h$, and $\theta_{T} \eqq \min\{  1, \theta_{T}^+ , \theta_{T}^- \}$, with 
    \begin{equation*}
        \theta_{T}^+ \eqq \begin{cases}
        \frac{b-\bar{v}_T}{M_T-\bar{v}_T} &\text{ if } M_T >b\\
        1 &\text{ otherwise}
    \end{cases}, \quad 
    \theta_{T}^- \eqq \begin{cases}
        \frac{\bar{v}_T-a}{\bar{v}_T-m_T} &\text{ if } m_T <a\\
        1 &\text{ otherwise}
    \end{cases}.
    \end{equation*}
    Then we set
    \begin{equation*}
        \Pi_h v |_T \eqq \bar{v}_T + \theta_{T} (\Tilde{v}_h-\bar{v}_T).
    \end{equation*}
    We show that $\Pi_h v \in X_h^{a,b,m}$. By construction, $\theta_{T}\in [0,1]$, hence
    \begin{equation*}
        \int_T \Pi_h v =\int_T (1-\theta_{T})\bar{v}_T + \theta_{T} \Tilde{v}_h = \int_T v,
    \end{equation*}
    which is the expected local mass-conservation which implies the global mass-conservation. Then, we recall the definition of $\theta_{T}$ to find
    \begin{equation*}
        a \leq \bar{v}_T + \theta_{T} (m_T-\bar{v}_T) \leq \Pi_h v|_T \leq \bar{v}_T + \theta_{T} (M_T-\bar{v}_T) \leq b,
    \end{equation*}
    which implies the nodal bound-preservation, hence $\Pi_h v \in X_h^{a,b,m}$.

    We then derive the convergence in different norms. For the $L^\infty$ norm, we want to show
    \begin{equation*}
        \norm{\Pi_h v - \Tilde{v}_h}_{L^\infty(T)} \leq C \norm{v-\Tilde{v}_h}_{L^\infty(T)},
    \end{equation*}
    which leads to 
    \begin{equation}\label{eq_opt_L_infty}
        \norm{v-\Pi_h v}_{L^\infty(T)} \leq \norm{v-\Tilde{v}_h }_{L^\infty(T)} + \norm{\Tilde{v}_h-\Pi_h v}_{L^\infty(T)} \leq C \norm{v-\Tilde{v}_h}_{L^\infty(T)}.
    \end{equation} 
    We distinguish the cases in which the minimum is attained by $\theta_{T}=1,\theta_{T}^+\neq 1,\theta_{T}^-\neq 1$, respectively. The case $\theta_{T}=1$ is trivial. If $\theta_{T}=\theta_{T}^+\neq 1$, we invoke Lemma \ref{lemma_poly_infty_min_max} and notice that $\int_T \Tilde{v}_h-\bar{v}_T =0$ to find
    \begin{align*}
        \norm{\Pi_h v - \Tilde{v}_h}_{L^\infty(T)} &= (1-\theta_{T}) \norm{\Tilde{v}_h-\bar{v}_T}_{L^\infty(T)} = (1-\theta_{T}^+) \norm{\Tilde{v}_h-\bar{v}_T}_{L^\infty(T)} \\
        &= \frac{M_T-b}{M_T-\bar{v}_T} \norm{\Tilde{v}_h-\bar{v}_T}_{L^\infty(T)} \leq C (M_T-b) \leq C \norm{v - \Tilde{v}_h}_{L^\infty(T)}.
    \end{align*}
    In the other case $\theta_{T}=\theta_{T}^-\neq 1$, similar argument leads to
    \begin{equation*}
        \norm{\Pi_h v - \Tilde{v}_h}_{L^\infty(T)} = (1-\theta_{T}^-) \norm{\Tilde{v}_h-\bar{v}_T}_{L^\infty(T)} \leq C\frac{|m_T-a|}{|\bar{v}_T-m_T|} |m_T-\bar{v}_T| \leq C \norm{v - \Tilde{v}_h}_{L^\infty(T)}.
    \end{equation*}
    Hence we complete the proof of the convergence in the $L^\infty$ norm, by using the approximation properties of $\Tilde{v}_h$.

    Now, we consider the convergence rate in the $H^k$-seminorms with $k=1,...,p+1$ and $v\in H^{p+1}(\calT_h)$.
    Repeated application of the discrete inverse inequality \eqref{eq_dis_trace_inv}, together with the $L^\infty$ stability of the $L^2$ projection and the Bramble–Hilbert estimate relating the $L^\infty$ error to the $H^{s+1}$ seminorm for $s+1>d/2$ in \cite[Proposition 4.3.2]{Brenner_Scott_2008}, yields
    \begin{equation*}
        \norm{\nabla^k (\Pi_h v - \Tilde{v}_h)}_T \leq C h_T^{d/2-k} \norm{ v - \Tilde{v}_h}_{L^\infty(T)} \leq C h_T^{s+1-k} \norm{\nabla^{s+1} v}_T.
    \end{equation*}
    Taking $s=p>\frac{d}{2}-1$ for $d\leq 3$ and $p\geq 1$, we conclude the proof of convergence for sufficiently smooth solutions.

    For the convergence in the $L^2$ norm, a direct calculation reveals $\norm{\Pi_h v - \Tilde{v}_h}_T \leq |T|^\frac{1}{2} \norm{\Pi_h v - \Tilde{v}_h}_{L^\infty(T)}$.
    This, with $ \norm{\Pi_h v - v}_T \leq  \norm{\Pi_h v - \Tilde{v}_h}_T +  \norm{\Tilde{v}_h-v}_T$, leads to the expected convergence under minimal regularity assumptions. For the convergence rate, we apply the Poincaré inequality with $\int_T v = \int_T\Pi_h v$ to find
    \begin{equation*}
        \norm{v-\Pi_h v}_T \leq C h_T \norm{\nabla (v-\Pi_h v)}_T \leq C h_T^{p+1} \norm{\nabla^{p+1} v}_T.
    \end{equation*}

    Under minimal regularity assumptions, the convergence in the $H^1$ seminorm needs a different argument. We denote the set of active cells by 
    \begin{equation*}
        \calA_h\eqq \{ T\in\calT_h \ : \ M_T>b \text{ or } m_T<a \},
    \end{equation*}
    and the corresponding active region by $A_h\eqq \bigcup_{T\in\calA_h}T$.
    We also define the modulus of continuity for $v$ by $\omega(\rho)\eqq \sup_{\xcoord,\boldsymbol{y}\in\bar{\Dom}, \ |\xcoord-\boldsymbol{y}|\leq \rho}  |v(\xcoord)-v(\boldsymbol{y})|$ for $\rho>0$, $\epsilon_h \eqq \max_{T\in\calT_h} \norm{v-\Tilde{v}_h}_{L^\infty(T)}$, and $\delta_h \eqq \epsilon_h + \omega(h)$. Since the $L^2$ projection is stable in the $L^\infty$ norm, we have $\epsilon_h\leq C \omega(h)\to 0$ as $h\to 0$, and consequently $\delta_h \to 0$. Next, we define the $\delta$-relaxed contact set by
    \begin{equation*}
        S_{\delta} \eqq \{ \xcoord\in \Dom \ : \ 0\leq v(\xcoord)-a \leq \delta \text{ or } 0\leq b-v(\xcoord) \leq \delta \}.
    \end{equation*}
    The above definition with decompositions $v(\xcoord)-a = (v(\xcoord) - m_T ) + (m_T-a)$ and $b-v(\xcoord) =  (b - M_T) + (M_T - v(\xcoord))$ for arbitrary $\xcoord\in \bar{T}$ with $T\in\calA_h$, implies $A_h\subset S_{\delta_h}$.
    Indeed, assume that $\xcoord\in A_h$ with $\xcoord \in T\in \calT_h$ and $m_T<a$, we have
    \begin{equation*}
        0\leq v(\xcoord)-a\leq v(\xcoord)-m_T +m_T-a\leq ( v(\xcoord)-v(\boldsymbol{y}) ) + ( v(\boldsymbol{y}) - \Tilde{v}_h(\boldsymbol{y}) ) \leq \omega(h) + \epsilon_h = \delta_h,
    \end{equation*}
    with $\boldsymbol{y}$ chosen such that $\Tilde{v}_h(\boldsymbol{y}) = m_T$. The same argument also holds for the case $M_T>b$.
    Then, the stability of the $L^2$ projection in the $H^1$ seminorm leads to
    \begin{equation*}
        \norm{\nabla_h(\Tilde{v}_h-\Pi_h v)}_\Dom^2 =\sum_{T\in\calA_h}(1-\theta_T)^2 \norm{\nabla \Tilde{v}_h}_T^2\leq C \int_{S_{\delta_h}}|\nabla v|^2. 
    \end{equation*}
    As $h\to 0$, we have $\delta_h\to 0$, hence $ \chi_{S_{\delta_h}} \to \chi_{\{ v=a \} \cup \{ v=b \} }$ a.e.. Since $v\in H^1(\Dom)$ satisfies $\nabla v=0$ a.e. on each level set (see, for example \cite[Lemma 7.7]{Gilbarg_Trudinger_1983}), the dominated convergence theorem implies
    \begin{equation*}
        \norm{\nabla_h (\Tilde{v}_h-\Pi_h v)}_\Dom^2 \leq C \int_{S_{\delta_h}}|\nabla v|^2 \to C \int_{\{ v=a \} \cup \{ v=b \}} |\nabla v|^2 =0 \text{ as } h\to 0.
    \end{equation*}
    This, with the triangle inequality leads to $\norm{\nabla_h (v-\Pi_h v)}_{\Dom} \to 0$ as $h\to 0$.

\end{proof}

\begin{remark}[Optimal candidate in $X_h^{a,b}$]
    We notice that, since $X_h^{a,b,m}\subset X_h^{a,b}$, the optimal candidate $\Pi_h v\in X_h^{a,b,m}$ is an admissible element in $X_h^{a,b}$, with the same optimal approximation properties. An alternative candidate in $X_h^{a,b}$ can be constructed by Lagrange interpolation.
\end{remark}

\subsection{A priori error estimate}
In this part, we consider the error estimates for the prediction step \eqref{eq_scheme_step1} and the correction step \eqref{eq_scheme_step2}. The error estimate of the prediction step can be derived by the standard argument, hence we only state the results and cite the proof given in the literature.
\begin{lemma}[Error estimate on $u-\bar{u}_h$]\label{lemma_dG_convergence}
    For the prediction step \eqref{eq_scheme_step1} with sufficiently large $\gamma_0>0$, we have
    \begin{equation}\label{eq_pre_H1_Cea}
        \dgnorm{u-\bar{u}_h} \leq C \inf_{v_h\in V_h^p} \Big[ \dgnorm{u-v_h} + \left( \sum_{T\in\calT_h} h_T\norm{K^\frac{1}{2}\nabla_h (u-v_h)}_{\partial T}^2 \right)^\frac{1}{2} \Big].
    \end{equation}
    Moreover, for $u\in H^{p+1}(\calT_h)$, we have
    \begin{equation}\label{eq_pre_H1_bound}
        \dgnorm{u-\bar{u}_h} \leq C \left(\sum_{T\in\calT_h} (\mu h_T^2 + K_{\max,T})h_T^{2p}\norm{\nabla^{p+1} u}_T^2\right)^\frac{1}{2}.
    \end{equation}
\end{lemma}

\begin{proof}
    The Céa-type result \eqref{eq_pre_H1_Cea} is proved in \cite[Theorem 3.1]{Ern_Stephansen_Zunino_2009}, and the optimal convergence result \eqref{eq_pre_H1_bound} is a direct consequence of Céa-type result \eqref{eq_pre_H1_Cea} and the continuous trace inequality \eqref{eq_con_trace}, as proved in \cite[Remark 3.1]{Ern_Stephansen_Zunino_2009}.
\end{proof}

\begin{lemma}[Error estimate on $u_h-\bar{u}_h$]\label{lemma_optimization_convergence}
    For the correction step \eqref{eq_scheme_step2}, we have
    \begin{equation}
        \dgnorm{u_h-\bar{u}_h} \leq C\left(\inf_{v_h\in X_h} \dgnorm{u-v_h} + \dgnorm{u-\bar{u}_h}\right).
    \end{equation}
\end{lemma}

\begin{proof}
    Since $u_h$ is the minimizer of the functional $a_h(\cdot-\bar{u}_h,\cdot-\bar{u}_h)$, we have
    \begin{equation*}
        C \dgnorm{u_h-\bar{u}_h}^2 \leq a_h(u_h-\bar{u}_h,u_h-\bar{u}_h) \leq a_h(v_h-\bar{u}_h,v_h-\bar{u}_h) \leq C \dgnorm{v_h-\bar{u}_h}^2 \quad \forall v_h\in X_h,
    \end{equation*}
    where we used the norm equivalence \eqref{eq_coercivity}. Then, the triangle inequality implies
    \begin{equation*}
        \dgnorm{u_h-\bar{u}_h}^2 \leq C \dgnorm{v_h-\bar{u}_h}^2 \leq C (\dgnorm{v_h-u}^2 + \dgnorm{u-\bar{u}_h}^2) \quad \forall v_h\in X_h,
    \end{equation*}
    which concludes the proof.
\end{proof}

\begin{theorem}[A priori error estimate]\label{thm_apriori_energy}
    Assume that $u\in H^{p+1}(\calT_h)$. Let $\gamma_0$ be sufficiently large and $X_h\in \{X_h^{a,b,m},X_h^{a,b}\}$. Then, for scheme \eqref{eq_scheme_step1}-\eqref{eq_scheme_step2}, we have the optimal convergence rate
    \begin{equation}
        \dgnorm{u-u_h}^2 \leq C \sum_{T\in\calT_h} h_T^{2p} (h_T^2\mu + K_{\max,T}) \norm{\nabla^{p+1}u}_{T}^2.
    \end{equation}
\end{theorem}

\begin{proof}
    The proof is a direct consequence of Lemmas \ref{lemma_opt_vh_abm}, \ref{lemma_dG_convergence} and \ref{lemma_optimization_convergence}. Indeed, the triangle inequality leads to
    \begin{align*}
        \dgnorm{u-u_h} \leq \dgnorm{u-\bar{u}_h} + \dgnorm{\bar{u}_h-u_h}.
    \end{align*}
    The first term on the right hand side is controlled by Lemma \ref{lemma_dG_convergence}, whereas the second term on the right hand side is controlled by  Lemmas \ref{lemma_opt_vh_abm}, \ref{lemma_dG_convergence} and \ref{lemma_optimization_convergence}.
\end{proof}

We finish this section with some remarks on minimal regularity and the extension to general linear problems.
\begin{remark}[Minimal regularity]
    The above theorem assumes $u\in H^2(\calT_h)$ to apply the standard analysis framework. However, such regularity assumption is only used to measure the error $\dgnorm{u-\bar{u}_h}$ of the prediction step \eqref{eq_scheme_step1}, and is independent of the correction step \eqref{eq_scheme_step2}. In \cite{Ern_Guermond_2022}, such regularity assumption is relaxed to $u\in H^{1+s}(\Dom)$ for arbitrarily small $s>0$, based on a face to cell lifting technique.
\end{remark}

\begin{remark}[General linear problems]
    We notice that the correction step only influences the error estimate on $\dgnorm{u_h-\bar{u}_h}$, and the error $\dgnorm{u-\bar{u}_h}$ is independent of the correction step. Hence, the a priori argument can be extended to any problem for which the prediction step admits an optimal estimate in a symmetric discrete energy norm and the admissible set contains a candidate with optimal approximation properties in that norm.
\end{remark}

\section{A posteriori analysis}\label{sec_aposteriori}

The key idea of our a posteriori analysis is to separate the prediction and correction errors using the triangle inequality. We estimate the prediction error by a standard energy argument and the correction error using the approximation properties of $\bar{u}_h$ and the triangle inequality.

Throughout this section, we assume for simplicity that
\begin{equation}\label{eq_aposteriori_data}
    g_D\in\calC(\bar{\Gamma}_D),\qquad g_D|_F\in\mathbb P_p(F)\quad\forall F\in\calF_D,\qquad
    g_N|_F\in\mathbb P_p(F)\quad\forall F\in\calF_N.
\end{equation}
The polynomial degree $p\geq1$ is fixed, and $\gamma_0$ is sufficiently large for \eqref{eq_coercivity}. Constants denoted by $C(\mu,K)$ are independent of the mesh size and may depend on $\mu$, $K_{\min}$ and $K_{\max}$, in addition to the dependencies allowed for $C$. We set
\begin{equation*}
    \chi_K\eqq\frac{K_{\max}}{K_{\min}}.
\end{equation*}
For $v\in H^1(\calT_h)$, we define
\begin{equation}\label{eq_local_dg_norm}
    \dgnorm{v}_T^2\eqq\mu\norm{v}_T^2+\norm{K^{1/2}\nabla_hv}_T^2
    +\frac12\sum_{F\in\calF_i\cap\partial T}\frac{\sigma_F}{h_F}\norm{\jump{v}}_F^2
    +\sum_{F\in\calF_D\cap\partial T}\frac{\sigma_F}{h_F}\norm{v}_F^2.
\end{equation}
Let $\calT_T\eqq\bigcup_{F\subset\partial T}\calT_F$ denote the set of cells sharing faces with given $T\in\calT_h$, and let $\Pi_0:L^2(\Dom)\to V_h^0$ denote the cellwise $L^2$ projection defined by
\begin{equation}\label{eq_cellmean_projection}
    (\Pi_0v)|_T\eqq\frac1{|T|}\int_Tv\quad\forall T\in\calT_h.
\end{equation}
The following parameters will be used later:
\begin{align}\label{eq_aposteriori_weights}
    \alpha_{1,T}&\eqq\min\left\{\frac{h_T}{\sqrt{K_{\min,T}}},\frac1{\sqrt\mu}\right\},&
    \alpha_{2,T}^2&\eqq\frac{K_{\max,T}+\mu h_T^2}{K_{\min,T}},\\
    \alpha_{3,T}&\eqq\min\left\{\sqrt{\frac{h_T}{K_{\min,T}}},\frac1{(\mu K_{\min,T})^{1/4}}\right\},&
    \alpha_{4,T}^2&\eqq\frac{K_{\max,T}}{h_T}\alpha_{3,T}^2.\nonumber
\end{align}
For $v_h\in V_h^p$, the local error estimators are
\begin{subequations}\label{eq_residual_indicators}
\begin{align}
    \eta_{R,T}(v_h)^2&\eqq\alpha_{1,T}^2\norm{(I-\Pi_0)\big(f+\DIV(K\nabla_hv_h)-\mu v_h\big)}_T^2,\\
    \eta_{J,T}(v_h)^2&\eqq\frac{\alpha_{3,T}^2}{2}\sum_{F\in\calF_i\cap\partial T}\norm{\jump{K\nabla_hv_h}\cdot\bn_F}_F^2,\\
    \eta_{N,T}(v_h)^2&\eqq\alpha_{3,T}^2\sum_{F\in\calF_N\cap\partial T}\norm{g_N-K\nabla_hv_h\cdot\bn_F}_F^2,\\
    \eta_{i,T}(v_h)^2&\eqq\frac{\alpha_{2,T}^2+\alpha_{4,T}^2}{2}\sum_{F\in\calF_i\cap\partial T}\frac{\sigma_F}{h_F}\norm{\jump{v_h}}_F^2,\\
    \eta_{D,T}(v_h)^2&\eqq(\alpha_{2,T}^2+\alpha_{4,T}^2)\sum_{F\in\calF_D\cap\partial T}\frac{\sigma_F}{h_F}\norm{v_h-g_D}_F^2.
\end{align}
\end{subequations}
We collect above estimators to define
\begin{align}
    &\eta_{0,T}^2\eqq\eta_{R,T}(\bar u_h)^2+\eta_{J,T}(\bar u_h)^2+\eta_{N,T}(\bar u_h)^2
        +\chi_K\big(\eta_{i,T}(\bar u_h)^2+\eta_{D,T}(\bar u_h)^2\big),\quad \eta_0^2\eqq\sum_{T\in\calT_h}\eta_{0,T}^2,\\
    &\eta_{dG,T}^2\eqq\eta_{0,T}^2+\dgnorm{u_h-\bar u_h}_T^2,\qquad
    \eta_{dG}^2\eqq\sum_{T\in\calT_h}\eta_{dG,T}^2
    =\eta_0^2+\dgnorm{u_h-\bar u_h}^2.\label{eq_eta_energy}
\end{align}
The data oscillation term is given by
\begin{equation}\label{eq_data_oscillation}
    \calR_{osc,T}^2\eqq\alpha_{1,T}^2\norm{f-\Pi_h^{L^2}f}_T^2,\qquad
    \calR_{osc}^2\eqq\sum_{T\in\calT_h}\calR_{osc,T}^2.
\end{equation}

The main theorem is given as follows, whose proof is postponed in the following sections.
\begin{theorem}[A posteriori error estimate with dG energy correction]\label{thm_aposteriori}
    Let $u_h$ be given by \eqref{eq_scheme_step1}--\eqref{eq_scheme_step2}. Under \eqref{eq_aposteriori_data}, we have
    \begin{equation}\label{eq_reliability_bound}
        \dgnorm{u-u_h}^2\leq C\eta_{dG}^2,
    \end{equation}
    and
    \begin{equation}\label{eq_efficiency_bound}
        \eta_{dG}^2\leq C(\mu,K)\big(\dgnorm{u-u_h}^2+\calR_{osc}^2\big).
    \end{equation}
    Moreover, for every $T\in\calT_h$,
    \begin{equation}\label{eq_local_efficiency_correction}
        \eta_{dG,T}^2\leq C(\mu,K)\sum_{T'\in\calT_T}
        \big(\dgnorm{u-u_h}_{T'}^2+\dgnorm{u_h-\bar u_h}_{T'}^2+\calR_{osc,T'}^2\big).
    \end{equation}
\end{theorem}

Then, we state a few remarks concerning the variants and practical implementation issues of the above theorem.

\begin{remark}[Nonconforming mesh]
    We assumed that the mesh $\calT_h$ is simplicial and conforming, hence there is no hanging node in our meshes. However, such assumption is used to simplify the analysis and is not mandatory. This assumption can be weakened to general mesh admitting a matching conforming shape-regular sub-triangulation. For details, we refer to \cite{DiPietro_Ern_2012}.
\end{remark}

\begin{remark}[Parameter robustness]
    For an augmented dG norm, it is possible to show the $K$-robustness, we refer to \cite{Ern_Stephansen_Voralik_2010} for details.
\end{remark}

\begin{remark}[Localized efficiency bound]
    We note that in the current result, the efficiency bound of error estimator $\eta_{dG,T}$ is not fully localized, with respect to the dG norm of $u-u_h$. This problem can be alleviated by using Lagrange multiplier $\sigma_h$ from \eqref{eq_sigma_h} to rewrite the prediction-correction scheme \eqref{eq_scheme_step1}-\eqref{eq_scheme_step2} as
    \begin{equation*}
        a_h(u_h,v_h) = l_h(v_h) + (\sigma_h, v_h) \quad \forall v_h\in V_h^p,
    \end{equation*}
    and then applying the variational inequality techniques (e.g., \cite{Braess_2005,Veeser_2001,Han_2024}) to get a localized efficiency bound. We do not pursue such localization here, owing to the numerical difficulty of evaluating $\sup_{v\in V_0,v\neq 0} \frac{\langle \sigma_h , v \rangle}{\norm{v}_{E,\Dom}}$, and the numerical findings in \S\ref{sec_num_test}. In \cite{Veeser_2001}, a computable residual estimator incorporating local terms of the form $h_T^2\norm{\nabla \sigma_h}_T$ is derived for a continuous piecewise affine discretization of the obstacle problem, with $\sigma_h$ defined through mass lumping, $K=I_d$, and homogeneous Dirichlet boundary conditions on $\partial\Dom=\Gamma_D$. However, extending this approach to mass-conservation constraints, anisotropic diffusion, or our high-order setting is nontrivial.
\end{remark}

\subsection{Error control of prediction step}
Following the standard energy argument, we have the reliable and efficient error bounds for the prediction step \eqref{eq_scheme_step1}.
\begin{lemma}[A posteriori error bound for prediction step]\label{lemma_reliability_sta}
    Under \eqref{eq_aposteriori_data}, the solution $\bar u_h$ of \eqref{eq_scheme_step1} satisfies
    \begin{equation}\label{eq_prediction_reliability}
        \dgnorm{u-\bar u_h}^2\leq C\eta_0^2,
    \end{equation}
    and
    \begin{equation}\label{eq_residual_efficiency_global}
        \eta_0^2\leq C(\mu,K)\big(\dgnorm{u-\bar u_h}^2+\calR_{osc}^2\big).
    \end{equation}
    Moreover, for every $T\in\calT_h$,
    \begin{equation}\label{eq_prediction_local_efficiency}
        \eta_{0,T}^2\leq C(\mu,K)\sum_{T'\in\calT_T}
        \big(\dgnorm{u-\bar u_h}_{T'}^2+\calR_{osc,T'}^2\big).
    \end{equation}
\end{lemma}

\begin{proof}
    This lemma can be demonstrated through the standard energy argument, as in \cite[Lemmas 3.1, 3.3--3.4 and Theorem 3.2]{Ern_Stephansen_2008}. Since the argument is standard, we only state the main steps and differences from their proof.

    \paragraph*{Reliability bound.} 
    We start by extending bilinear form $a$ to $H^1(\calT_h)$ and set
    \begin{align}\label{eq_energy_extension}
        a(v,w)&\eqq\sum_{T\in\calT_h}\big((\mu v,w)_T+(K\nabla_hv,\nabla_hw)_T\big),\\
        \norm{v}_{E,S}^2&\eqq\sum_{T\in\calT_h}\int_{S\cap T}
        \big(\mu v^2+(K\nabla_hv)\cdot\nabla_hv\big).\nonumber
    \end{align}
    Let $\calI_{Os}:V_h^p\to V_h^p\cap V$ be the Oswald averaging operator with Dirichlet nodal values on $\Gamma_D$:
    \begin{equation*}
        \calI_{Os} v_h (s) = \begin{cases}
            \frac{1}{|\calT_s|}\sum_{T\in\calT_s} v_h |_T(s) &\text{ if } s\in\calV^{cG}, \ s\not\in \Bar{\Gamma}_D\\
            g_D(s)  &\text{ otherwise }
        \end{cases} \quad \forall v_h \in V_h^p,
    \end{equation*}
    with $\calT_s$ the set of mesh elements that contain the node $s$. Then, we recall the error control for the Oswald averaging operator (see, for example, \cite[(3.23)--(3.24)]{Ern_Stephansen_2008})
    \begin{equation}\label{eq_oswald_approx}
        \norm{\nabla^j(v_h-\calI_{Os} v_h)}_T^2 \leq C \Big( \sum_{F\in \Tilde{\calF}_T \cap \calF_i } h_F^{1-2j} \norm{\jump{v_h}}_F^2 + \sum_{F\in \Tilde{\calF}_T \cap \calF_D} h_F^{1-2j} \norm{v_h -g_D}_F^2 \Big) \quad j=0,1,
    \end{equation}
    and the fact that $\sigma_F\geq K_{\min}/2$ and $\alpha_{2,T}^2\geq1+\mu h_T^2/K_{\max}$, to obtain
    \begin{align}
        \norm{\bar u_h-\calI_{Os}\bar u_h}_{E,\Dom}^2
        &\leq C\sum_{F\in\calF_i}
        \left(\frac{K_{\max}}{h_F}+\mu h_F\right)
        \norm{\jump{\bar u_h}}_F^2\nonumber + C\sum_{F\in\calF_D}
        \left(\frac{K_{\max}}{h_F}+\mu h_F\right)
        \norm{\bar u_h-g_D}_F^2\nonumber\\
        &\leq C\chi_K\sum_{T\in\calT_h}
        \big(\eta_{i,T}(\bar u_h)^2+\eta_{D,T}(\bar u_h)^2\big).
        \label{eq_nonconforming_bound}
    \end{align}
    With the help of $w = u-\calI_{Os}\bar u_h\in V_0$, we find 
    \begin{align*}
        \norm{\bar u_h-\calI_{Os}\bar u_h}_{E,\Dom}^2
        & = \norm{(u-\bar u_h)-w}_{E,\Dom}^2\nonumber=\norm{u-\bar u_h}_{E,\Dom}^2
        -2a(u-\bar u_h,w)+\norm{w}_{E,\Dom}^2 \\
        & = \norm{u-\bar u_h}_{E,\Dom}^2
        -2( l(w)-a(\bar u_h,w) )+\norm{w}_{E,\Dom}^2,
    \end{align*}
    which leads to
    \begin{align}
        \norm{u-\bar u_h}_{E,\Dom}^2
        &=\norm{\bar u_h-\calI_{Os}\bar u_h}_{E,\Dom}^2
        +2\big(l(w)-a(\bar u_h,w)\big)
        -\norm{w}_{E,\Dom}^2\nonumber\\
        &\leq\norm{\bar u_h-\calI_{Os}\bar u_h}_{E,\Dom}^2
        +\left(
        \sup_{0\neq\varphi\in V_0}
        \frac{|l(\varphi)-a(\bar u_h,\varphi)|}
        {\norm{\varphi}_{E,\Dom}}
        \right)^2.
        \label{eq_abstract_reliability}
    \end{align}

    For $\varphi\in V_0$, set $w=\varphi-\Pi_0\varphi$. Testing \eqref{eq_scheme_step1} with $\Pi_0\varphi\in V_h^p$ and integrating by parts, we obtain
    \begin{align}\label{eq_conforming_residual}
        l(\varphi)-a(\bar u_h,\varphi)
        ={}&\sum_T\big((I-\Pi_0)(f+\DIV(K\nabla_h\bar u_h)-\mu\bar u_h),w\big)_T\\
        &-\sum_{F\in\calF_i}(\jump{K\nabla_h\bar u_h}\cdot\bn_F,\avgbar{w})_F
        +\sum_{F\in\calF_N}(g_N-K\nabla_h\bar u_h\cdot\bn_F,w)_F\nonumber\\
        &-\sum_{F\in\calF_i}\gamma_0\frac{\sigma_F}{h_F}(\jump{\bar u_h},\jump{w})_F
        -\sum_{F\in\calF_D}\gamma_0\frac{\sigma_F}{h_F}(\bar u_h-g_D,w)_F.\nonumber
    \end{align}
    where we used the orthogonality property
    \begin{equation}\label{eq_residual_mean_orthogonality}
        \big(f+\DIV(K\nabla_h\bar u_h)-\mu\bar u_h,w\big)_T
        =\big((I-\Pi_0)(f+\DIV(K\nabla_h\bar u_h)-\mu\bar u_h), w\big)_T.
    \end{equation}
    Then, the Poincaré inequality with the continuous trace inequality \eqref{eq_con_trace} (or estimates in \cite[(3.36), (3.38)]{Ern_Stephansen_2008}) yield
    \begin{equation}\label{eq_prediction_test_bounds}
        \norm{\varphi - \Pi_0\varphi}_T\leq C\alpha_{1,T}\norm{\varphi}_{E,T},\qquad
        \norm{\varphi - \Pi_0\varphi}_{\partial T}\leq C\alpha_{3,T}\norm{\varphi}_{E,T}.
    \end{equation}
    Using the above estimates, $\sigma_F\leq K_{\max,T}$, $\alpha_{3,T}^2\left(\gamma_0\frac{\sigma_F}{h_F}\right)^2
        \leq C\gamma_0^2\alpha_{4,T}^2\frac{\sigma_F}{h_F}$, 
    and the Cauchy--Schwarz inequality in \eqref{eq_conforming_residual}, we obtain
    \begin{equation}\label{eq_conforming_residual_bound}
        |l(\varphi)-a(\bar u_h,\varphi)|\leq C\eta_0\norm{\varphi}_{E,\Dom}.
    \end{equation}
    Inserting \eqref{eq_nonconforming_bound} and \eqref{eq_conforming_residual_bound} in \eqref{eq_abstract_reliability}, and adding the jump contribution proves \eqref{eq_prediction_reliability}.

    \paragraph*{Efficiency bound.} 
    For efficiency, we follow the standard element and face bubble arguments of \cite[Propositions 3.3--3.4]{Ern_Stephansen_2008}. For any $v_h\in V_h^p$, let $b_T$ be an element bubble (see, for example \cite{Verfurth_2013} and \cite[Lemma 5.27]{DiPietro_Ern_2012}) and set
    \begin{equation*}
        \psi_T\eqq b_T\big(\Pi_h^{L^2}f+\DIV(K\nabla_hv_h)-\mu v_h\big)
        \in H_0^1(T).
    \end{equation*}
    which satisfies
    \begin{align*}
        \norm{\psi_T}_T^2 +h_T^2\norm{\nabla\psi_T}_T^2 \leq C \norm{\Pi_h^{L^2}f+\DIV(K\nabla_hv_h)-\mu v_h}_T^2
        \leq C\big(\Pi_h^{L^2}f+\DIV(K\nabla_hv_h)-\mu v_h,\psi_T\big)_T.
    \end{align*}
    Then, integration by parts gives
    \begin{align*}
        &\quad \big(\Pi_h^{L^2}f+\DIV(K\nabla_hv_h)-\mu v_h,\psi_T\big)_T
        =a(u-v_h,\psi_T)+(\Pi_h^{L^2}f-f,\psi_T)_T\\
        &\leq \norm{u-v_h}_{E,T}
        \left(\sqrt{K_{\max,T}}\norm{\nabla\psi_T}_T
        +\sqrt\mu\norm{\psi_T}_T\right) + \norm{f-\Pi_h^{L^2}f}_T\norm{\psi_T}_T.
    \end{align*}
    Consequently,
    \begin{align*}
        &\norm{\Pi_h^{L^2}f+\DIV(K\nabla_hv_h)-\mu v_h}_T \leq C\left[
        \left(\frac{\sqrt{K_{\max,T}}}{h_T}+\sqrt\mu\right)
        \norm{u-v_h}_{E,T}+\norm{f-\Pi_h^{L^2}f}_T
        \right].
    \end{align*}
    Since $I-\Pi_0$ is an $L^2$ contraction and $\alpha_{1,T}
        \left(\frac{\sqrt{K_{\max,T}}}{h_T}+\sqrt\mu\right)
        \leq 2\sqrt{\frac{K_{\max,T}}{K_{\min,T}}}$, 
    the triangle inequality yields
    \begin{align}
        \eta_{R,T}(v_h)^2 \leq\alpha_{1,T}^2
        \norm{f+\DIV(K\nabla_hv_h)-\mu v_h}_T^2 \leq C\frac{K_{\max,T}}{K_{\min,T}}
        \norm{u-v_h}_{E,T}^2+C\calR_{osc,T}^2.
        \label{eq_efficiency_residual}
    \end{align}
    
    For $F\in\calF_i$, let
    $T_F\eqq\operatorname{int}(\bigcup_{T'\in\calT_F}\bar T')$. The standard face bubble function is given by $b_F$, and we define $r_F = \jump{K\nabla_hv_h}\cdot\bn_F$, extend $r_F$ to $T_F$ by constant values along $\bn_F$, and set $\psi_F=b_F r_F\in H_0^1(T_F)$. Such construction satisfies
    \begin{align*}
        \norm{\jump{K\nabla_hv_h}\cdot\bn_F}_F^2
        \leq C(\jump{K\nabla_hv_h}\cdot\bn_F,\psi_F)_F, \quad
        \norm{\psi_F}_{T_F}
        +h_F\norm{\nabla\psi_F}_{T_F}
        \leq Ch_F^{1/2}
        \norm{\jump{K\nabla_hv_h}\cdot\bn_F}_F.
    \end{align*}
    Then, cellwise integration by parts implies
    \begin{equation*}
        (\jump{K\nabla_h v_h}\cdot\bn_F,\psi_F)_F
        =\sum_{T'\in\calT_F}
        \big(f+\DIV(K\nabla_hv_h)-\mu v_h,\psi_F\big)_{T'}
        -a(u-v_h,\psi_F).
    \end{equation*}
    This leads to
    \begin{align}
        &\quad C\norm{\jump{K\nabla_hv_h}\cdot\bn_F}_F\nonumber\\
        &\leq h_F^{1/2}\sum_{T'\in\calT_F}
        \norm{f+\DIV(K\nabla_hv_h)-\mu v_h}_{T'} + \sum_{T'\in\calT_F}
        \left(\sqrt{K_{\max,T'}}h_F^{-1/2}
        +\sqrt\mu h_F^{1/2}\right)
        \norm{u-v_h}_{E,T'}.
        \label{eq_face_bubble_bound}
    \end{align}
    For $F\subset\partial T$ and $T'\in\calT_F$, shape regularity gives
    $h_T\simeq h_F\simeq h_{T'}$. The definitions of the weights imply
    \begin{equation*}
        \alpha_{3,T}\left(
        \frac{h_F^{1/2}}{\alpha_{1,T'}}
        +\sqrt{K_{\max,T'}}h_F^{-1/2}
        +\sqrt\mu h_F^{1/2}\right)
        \leq C\left(
        \sqrt{\chi_K}
        +\sqrt{\frac{\mu}{K_{\min}}}\diam(\Dom)\right).
    \end{equation*}
    Multiplying \eqref{eq_face_bubble_bound} by $\alpha_{3,T}$,
    using \eqref{eq_efficiency_residual}, and summing over the interior faces of $T$, we find
    \begin{equation}\label{eq_efficiency_flux_jump}
        \eta_{J,T}(v_h)^2
        \leq C(\mu,K)\sum_{T'\in\calT_T}
        \big(\norm{u-v_h}_{E,T'}^2+\calR_{osc,T'}^2\big).
    \end{equation}
    The estimate for $\eta_{N,T}(v_h)$ follows similarly, hence we omit for brevity. Since the jump estimators are included in the dG energy norm, we have
    \begin{align}\label{eq_residual_efficiency_arbitrary}
        &\alpha_{1,T}^2\norm{f+\DIV(K\nabla_hv_h)-\mu v_h}_T^2
        +\eta_{J,T}(v_h)^2+\eta_{N,T}(v_h)^2\\
        &\qquad+\eta_{i,T}(v_h)^2+\eta_{D,T}(v_h)^2
        \leq C(\mu,K)\sum_{T'\in\calT_T}
        \big(\dgnorm{u-v_h}_{T'}^2+\calR_{osc,T'}^2\big).\nonumber
    \end{align}
    Finally, we take $v_h=\bar u_h$ to conclude \eqref{eq_residual_efficiency_global}.
\end{proof}

\subsection{Error control of the prediction-correction scheme}
For the error control of the correction step, we consider the triangle inequality $\dgnorm{u-u_h}\leq \dgnorm{u-\bar{u}_h} + \dgnorm{\bar{u}_h-u_h}$ and the fact that $\dgnorm{\bar{u}_h-u_h}$ and $\dgnorm{\bar{u}_h-u}$ are sufficiently close, up to a data oscillation term. See, for examples \cite[Lemma 2.1 and Section 3.2]{Gudi_2010} and \cite[Lemma 5.27]{DiPietro_Ern_2012}.

\begin{lemma}[Best approximation of prediction step]\label{lemma_prediction_best_approximation}
    For every $v_h\in V_h^p$, we have
    \begin{equation}\label{eq_prediction_comparison}
        \dgnorm{\bar u_h-v_h}
        \leq C(\mu,K)\big(\dgnorm{u-v_h}+\calR_{osc}\big).
    \end{equation}
    Consequently,
    \begin{equation}\label{eq_prediction_best_approximation}
        \dgnorm{u-\bar u_h}\leq C(\mu,K)
        \left(\inf_{v_h\in V_h^p}\dgnorm{u-v_h}+\calR_{osc}\right).
    \end{equation}
\end{lemma}

\begin{proof}
    Let $w_h=\bar u_h-v_h$, $z_h=w_h-\calI_{Os}^0w_h$, where $\calI_{Os}^0$ is the Oswald operator used in the proof of Lemma \ref{lemma_reliability_sta}, but takes zero value at Dirichlet boundary nodes. The coercivity of $a_h$ \eqref{eq_coercivity} and the prediction scheme \eqref{eq_scheme_step1} now imply
    \begin{align*}
        C\dgnorm{w_h}^2
        &\leq a_h(w_h,w_h)=l_h(w_h)-a_h(v_h,w_h)\\
        &= \left[l_h(\calI_{Os}^0w_h)-a_h(v_h,\calI_{Os}^0w_h) \right] + \left[ l_h(z_h)-a_h(v_h,z_h) \right],
    \end{align*}
    which will be controlled separately.
    
    The $L^2$ stability of $\calI_{Os}^0$ and the approximation properties \eqref{eq_oswald_approx} imply
    \begin{equation}\label{eq_oswald_global}
        \norm{\calI_{Os}^0w_h}_{E,\Dom}^2+\dgnorm{z_h}^2\leq C\chi_K\dgnorm{w_h}^2,
    \end{equation}
    and
    \begin{equation}\label{eq_oswald_scaled}
        \sum_T\left(h_T^{-2}\norm{z_h}_T^2+h_T^{-1}\norm{z_h}_{\partial T}^2\right)
        \leq\frac{C}{K_{\min}}\dgnorm{w_h}^2.
    \end{equation}
    In particular, the definitions of $\alpha_{1,T}$ and $\alpha_{3,T}$ imply
    \begin{equation}\label{eq_discrete_residual_test_bounds}
        \sum_T\left(\alpha_{1,T}^{-2}\norm{z_h}_T^2+
        \alpha_{3,T}^{-2}\norm{z_h}_{\partial T}^2\right)
        \leq C(\mu,K)\dgnorm{w_h}^2.
    \end{equation}
    Since $\calI_{Os}^0w_h\in V_0$, the weak formulation \eqref{eq_model_pb_weak} gives
    \begin{align*}
        l_h(\calI_{Os}^0w_h)-a_h(v_h,\calI_{Os}^0w_h)
        &= a(u-v_h,\calI_{Os}^0w_h) \\
        &\quad +\sum_{F\in\calF_i}(\jump{v_h},\avg{K\nabla \calI_{Os}^0w_h}\cdot\bn_F)_F +\sum_{F\in\calF_D}(v_h-g_D,K\nabla \calI_{Os}^0w_h\cdot\bn_F)_F.
    \end{align*}
    On an interior face, the harmonic weights satisfy $(\omega_F^+)^2\delta_F^++(\omega_F^-)^2\delta_F^-=\sigma_F$. Hence, the discrete trace inequality \eqref{eq_dis_trace_inv} yields
    \begin{equation}\label{eq_weighted_flux_trace_discrete}
        \sum_{F\in\calF_i\cup\calF_D}\frac{h_F}{\sigma_F}
        \norm{\avg{K\nabla_h\xi_h}\cdot\bn_F}_F^2
        \leq C\norm{K^{1/2}\nabla_h\xi_h}_\Dom^2
        \quad\forall\xi_h\in V_h^p.
    \end{equation}
    Using above estimate and \eqref{eq_oswald_global}, we find
    \begin{equation}\label{eq_discrete_residual_conforming}
        |l_h(\calI_{Os}^0w_h)-a_h(v_h,\calI_{Os}^0w_h)|
        \leq C(K)\dgnorm{u-v_h}\dgnorm{w_h}.
    \end{equation}

    For the remaining part, integration by parts yields
    \begin{align}\label{eq_discrete_residual_identity}
        l_h(z_h)-a_h(v_h,z_h)
        ={}&\sum_T\big(f+\DIV(K\nabla_hv_h)-\mu v_h,z_h\big)_T\\
        &-\sum_{F\in\calF_i}(\jump{K\nabla_hv_h}\cdot\bn_F,\avgbar{z_h})_F
        +\sum_{F\in\calF_N}(g_N-K\nabla_hv_h\cdot\bn_F,z_h)_F\nonumber\\
        &+\sum_{F\in\calF_i}\left(\jump{v_h},\avg{K\nabla_hz_h}\cdot\bn_F
        -\gamma_0\frac{\sigma_F}{h_F}\jump{z_h}\right)_F\nonumber\\
        &+\sum_{F\in\calF_D}\left(v_h-g_D,K\nabla_hz_h\cdot\bn_F
        -\gamma_0\frac{\sigma_F}{h_F}z_h\right)_F.\nonumber
    \end{align}
    Since $\alpha_{2,T}^2+\alpha_{4,T}^2\geq1$, the Cauchy--Schwarz inequality, \eqref{eq_discrete_residual_test_bounds} and \eqref{eq_weighted_flux_trace_discrete} imply
    \begin{align}\label{eq_discrete_residual_bound}
        |l_h(z_h)-a_h(v_h,z_h)|
        &\leq C(\mu,K)\Bigg(\sum_T\Big[
        \alpha_{1,T}^2\norm{f+\DIV(K\nabla_hv_h)-\mu v_h}_T^2\\
        &\qquad+\eta_{J,T}(v_h)^2+\eta_{N,T}(v_h)^2
        +\eta_{i,T}(v_h)^2+\eta_{D,T}(v_h)^2\Big]\Bigg)^{1/2}\dgnorm{w_h}\nonumber\\
        &\leq C(\mu,K)\big(\dgnorm{u-v_h}+\calR_{osc}\big)\dgnorm{w_h},\nonumber
    \end{align}
    where the last inequality follows from \eqref{eq_residual_efficiency_arbitrary}.

    This proves \eqref{eq_prediction_comparison}. \eqref{eq_prediction_best_approximation} is a direct consequence of the triangle inequality and the fact that $v_h\in V_h^p$ is arbitrary.
\end{proof}

Now, we are ready to combine the preceding results to prove the main theorem.
\begin{proof}[Proof of Theorem \ref{thm_aposteriori}]
    The triangle inequality and \eqref{eq_prediction_reliability} give
    \begin{equation*}
        \dgnorm{u-u_h}\leq\dgnorm{u-\bar u_h}+\dgnorm{u_h-\bar u_h}
        \leq C\eta_0+\dgnorm{u_h-\bar u_h}.
    \end{equation*}
    which proves the reliability bound \eqref{eq_reliability_bound}.

    Then we consider the efficiency bound. Taking $v_h=u_h$ in \eqref{eq_prediction_comparison}, and using the triangle inequality, yields
    \begin{align}\label{eq_correction_efficiency}
        \dgnorm{u_h-\bar u_h}&\leq C(\mu,K)\big(\dgnorm{u-u_h}+\calR_{osc}\big),\\
        \dgnorm{u-\bar u_h}&\leq C(\mu,K)\big(\dgnorm{u-u_h}+\calR_{osc}\big).\nonumber
    \end{align}
    Moreover, the efficiency bound of the prediction step \eqref{eq_residual_efficiency_global} implies
    \begin{equation*}
        \eta_0^2\leq C(\mu,K)\big(\dgnorm{u-u_h}^2+\calR_{osc}^2\big),
    \end{equation*}
    which leads to \eqref{eq_efficiency_bound}. Finally,  efficiency \eqref{eq_prediction_local_efficiency} and
    \begin{equation*}
        \dgnorm{u-\bar u_h}_{T'}^2
        \leq2\dgnorm{u-u_h}_{T'}^2+2\dgnorm{u_h-\bar u_h}_{T'}^2
    \end{equation*}
    give \eqref{eq_local_efficiency_correction}, which completes the proof.
\end{proof}

Next corollary is a direct application of the proof of Theorem \ref{thm_aposteriori} with $u_h$ replaced by $v_h$, and will be used to derive the a posteriori estimate in \S\ref{sec_L2_correction}:
\begin{corollary}[A posteriori error estimate with arbitrary numerical solution]\label{cor_arbitrary_correction}
    For every $v_h\in V_h^p$, define $\eta(v_h)^2\eqq\eta_0^2+\dgnorm{v_h-\bar u_h}^2$. Then
    \begin{equation}\label{eq_arbitrary_correction_bound}
        \dgnorm{u-v_h}^2\leq C\eta(v_h)^2,\qquad
        \eta(v_h)^2\leq C(\mu,K)\big(\dgnorm{u-v_h}^2+\calR_{osc}^2\big).
    \end{equation}
\end{corollary}

\section{Analysis of $L^2$ based scheme}\label{sec_L2_correction}
In the following sections, we establish the a priori and a posteriori error estimates of the $L^2$ correction based scheme \eqref{eq_scheme_step1}-\eqref{eq_scheme_step2_L2}.

\subsection{A priori analysis}

As in the standard argument, we need to apply the Aubin--Nitsche argument to derive the optimal $L^2$ error estimate, hence we consider the dual problem
\begin{equation}\label{eq_model_z}
    \mu z -\DIV(K\nabla z)= u-\bar{u}_h \text{ in } \Dom, \qquad z=0\text{ on }\Gamma_D,\qquad K\nabla z \cdot \bn_{\Dom} = 0 \text{ on }\Gamma_N.
\end{equation}
In addition, we assume that the broken $H^2$-regularity holds true for the above dual problem:
\begin{equation}\label{eq_H2_z}
    \norm{\nabla_h^2 z}_{\Dom} \leq C K_{\min}^{-1} \norm{u-\bar{u}_h}_\Dom,
\end{equation}
as considered in \cite{Ern_Stephansen_Zunino_2009}. With these tools, we are ready to state the convergence in $L^2$  and energy norms for the scheme \eqref{eq_scheme_step1}-\eqref{eq_scheme_step2_L2}:
\begin{theorem}[A priori estimate for $L^2$ correction]\label{Thm_apriori_L2}
    Let $\bar{u}_h$ and $\mathfrak{u}_h$ be given by \eqref{eq_scheme_step1} and \eqref{eq_scheme_step2_L2}, respectively. Assume that $u\in\calC(\bar{\Dom};[a,b])\cap H^{p+1}(\calT_h)$ and the broken $H^2$-regularity \eqref{eq_H2_z} holds true. Let $\gamma_0$ be sufficiently large and $X_h\in \{X_h^{a,b,m},X_h^{a,b}\}$. We define $\beta_T\eqq K_{\max,T}+\mu h_T^2$ and $\beta\eqq \max_{T\in\calT_h} \beta_T$. Then, we have
    \begin{align}
        \dgnorm{u-\mathfrak{u}_h}^2 &\leq C \sum_{T\in\calT_h} h_T^{2p} \Big( \beta_T + \frac{h_T^2}{h_{\min}^2}K_{\max} + \frac{h^2}{h_{\min}^2}\frac{K_{\max}\beta\beta_T}{K_{\min}^2} \Big) \norm{\nabla^{p+1} u}_T^2,\label{eq_error_L2_cor_H1}\\
        \norm{u-\mathfrak{u}_h}_\Dom^2 &\leq C \sum_{T\in\calT_h} h_T^{2p+2}  \left(1 +  \frac{h^2}{h_T^2}\frac{\beta \beta_T}{K_{\min}^2}\right)\norm{\nabla^{p+1} u}_T^2 .\label{eq_error_L2_cor_L2}
    \end{align}
\end{theorem}

From the above theorem, we can deduce the following results:
\begin{corollary}
    Under the assumptions of Theorem \ref{Thm_apriori_L2}, the structure-preserving numerical solution $\mathfrak{u}_h\in X_h$ from \eqref{eq_scheme_step1}-\eqref{eq_scheme_step2_L2} has the optimal convergence rate in  $L^2$ norm
    \begin{align}
        \norm{u-\mathfrak{u}_h} \leq C(u,K,\mu) h^{p+1}.
        \end{align}
    In addition, if $\calT_h$ is quasi-uniform in the sense that $h\leq C h_{\min}$, it also has the optimal convergence rate in the dG energy norm:
    \begin{align}
     \dgnorm{u-\mathfrak{u}_h} \leq C(u,K,\mu) h^{p}.
    \end{align}
    In the above $C(u,K,\mu)$ is a generic constant depending on $u$, $\mu$ and $K$.
\end{corollary}

\begin{remark}[$L^2$ error estimate in the literature]
    To the best of our knowledge, in the literature, the optimal $L^2$ convergence rate with respect to $h$ of bound-preserving schemes is proved only in two categories, for $p = 1$. The first one focuses on standard linear schemes under the assumption that a discrete maximum principle holds true (\eg \cite{Barrenechea_John_Knobloch_2024}), hence standard Aubin--Nitsche argument leads to the optimal convergence rate in the $L^2$ norm; the second one \cite{Barrenechea_ChaumontFrelet_2026}, instead, proposes an adaptation of the standard Aubin--Nitsche argument for the nonlinear scheme from \cite{Barrenechea_Georgoulis_Pryer_Veeser_2024}. Our methodology does not follow any of these, but relies only on the $L^2$ optimal error estimate of the prediction step \eqref{eq_scheme_step1} and the optimal approximation properties of the correction step \eqref{eq_scheme_step2_L2}. The price paid for this simple argument is the use of the discrete inverse inequality for the error estimate in the dG energy norm.
\end{remark}

\begin{proof}[Proof of Theorem \ref{Thm_apriori_L2}]
    We split the proof into three parts. In the first part, we show the optimal $L^2$ error estimate  for $\bar{u}_h$ by using Aubin--Nitsche argument; in the second part, we derive the optimal $L^2$ error estimate for $\mathfrak{u}_h$, by using the $L^2$ optimality condition of the correction step \eqref{eq_scheme_step2_L2}; in the last step, we recover the optimal convergence in the dG energy norm, with the help of the continuous and discrete inverse inequalities.

    \paragraph{Error estimate on $u-\bar{u}_h$ in the $L^2$ norm.} The proof essentially follows the same argument in \cite[Theorem 3.2, Corollary 3.1]{Ern_Stephansen_Zunino_2009}, hence we only give a sketch of the proof here. 
    Let $e=u-\bar{u}_h$ and let $z$ solve \eqref{eq_model_z}; since $z=0$ on $\Gamma_D$, $K\nabla z\cdot\bn_{\Dom}=0$ on $\Gamma_N$, and its normal flux is continuous across interior faces, adjoint consistency and Galerkin orthogonality give, for every $z_h\in V_h^p$,
    \begin{equation*}
        \norm{e}_{\Dom}^2=a_h(e,z)=a_h(e,z-z_h).
    \end{equation*}
    Picking $z_h$ as the conforming Lagrange interpolant of $z$, we have
    \begin{equation*}
        \norm{e}_{\Dom}^2
        \leq C\dgnorm{e}
        \left(
            \dgnorm{z-z_h}^2
            +\sum_{T\in\calT_h}h_T
            \norm{K^{1/2}\nabla(z-z_h)}_{\partial T}^2
        \right)^{1/2}.
    \end{equation*}
    The approximation properties of interpolation $z_h$ and trace estimates \eqref{eq_con_trace}, together with \eqref{eq_H2_z}, yield
    \begin{equation*}
        \dgnorm{z-z_h}^2
        +\sum_{T\in\calT_h}h_T
            \norm{K^{1/2}\nabla(z-z_h)}_{\partial T}^2
        \leq C h^2\frac{K_{\max}+\mu h^2}{K_{\min}^2}
            \norm{e}_{\Dom}^2.
    \end{equation*}
    The above estimate with \eqref{eq_pre_H1_bound} leads to
    \begin{equation}\label{eq_error_L2_pre}
        \norm{u-\bar{u}_h}_{\Dom}^2
        \leq C h^2 \frac{K_{\max}+\mu h^2}{K_{\min}^2}
        \sum_{T\in\calT_h} h_T^{2p}(K_{\max,T}+\mu h_T^2)
        \norm{\nabla^{p+1}u}_T^2
        \leq C\frac{h^2\beta}{K_{\min}^2}
        \sum_{T\in\calT_h}h_T^{2p}\beta_T
        \norm{\nabla^{p+1}u}_T^2.
    \end{equation}
    
    \paragraph{Error estimate on $u-\mathfrak{u}_h$ in the $L^2$ norm.} 
    Since $\mathfrak{u}_h$ is the minimizer of \eqref{eq_scheme_step2_L2}, the triangle inequality implies that for all $v_h\in X_h$,
    \begin{align*}
        \norm{u-\mathfrak{u}_h}_\Dom &\leq \norm{u-\bar{u}_h}_\Dom + \norm{\bar{u}_h-\mathfrak{u}_h}_\Dom \leq   \norm{u-\bar{u}_h}_\Dom + \norm{\bar{u}_h-v_h}_\Dom \leq 2\norm{u-\bar{u}_h}_\Dom + \norm{u-v_h}_\Dom.
    \end{align*}
    Setting $v_h= \Pi_h u$ defined in Lemma \ref{lemma_opt_vh_abm}, the approximation properties of $\Pi_h u$ yield
    \begin{align*}
        \norm{v_h-u}_\Dom^2 \leq C\sum_{T\in\calT_h} h_T^{2p+2} \norm{\nabla^{p+1} u}_T^2,
    \end{align*}
    which leads to the expected $L^2$ error bound \eqref{eq_error_L2_cor_L2}.

    \paragraph{Error estimate in the dG norm.}
    The error in the dG norm can be decomposed as 
    \begin{equation*}
        \dgnorm{u-\mathfrak{u}_h}^2 = \mu \norm{u-\mathfrak{u}_h}_\Dom^2 + \norm{K^\frac{1}{2}\nabla_h(u-\mathfrak{u}_h)}^2 + \sum_{F\in\calF_i\cup\calF_D} \frac{\sigma_F}{h_F}\norm{\jump{u-\mathfrak{u}_h}}_F^2.
    \end{equation*}
    We use the triangle inequality, Lemma \ref{lemma_dG_convergence}, discrete trace and inverse inequalities \eqref{eq_dis_trace_inv} and the fact that $\mathfrak{u}_h$ is the minimizer of \eqref{eq_scheme_step2_L2}, to derive, for $v_h=\Pi_h u \in X_h$,
    \begin{align*}
        \norm{\mu^\frac{1}{2}(u-\mathfrak{u}_h)}_{\Dom} &\leq \norm{\mu^\frac{1}{2}(u-\bar{u}_h)}_{\Dom} + \norm{\mu^\frac{1}{2}(\bar{u}_h-\mathfrak{u}_h)}_{\Dom} \leq 2 \norm{\mu^\frac{1}{2}(u-\bar{u}_h)}_{\Dom} +\norm{\mu^\frac{1}{2}(u-v_h)}_{\Dom},
    \end{align*}
    and $\norm{K^\frac{1}{2}\nabla_h(u-\mathfrak{u}_h)}_\Dom \leq \norm{K^\frac{1}{2}\nabla_h(u-\bar{u}_h)}_\Dom + \norm{K^\frac{1}{2}\nabla_h(\bar{u}_h-\mathfrak{u}_h)}_\Dom$, with 
    \begin{align*}
        \norm{K^\frac{1}{2}\nabla_h(\bar{u}_h-\mathfrak{u}_h)}_\Dom \leq C \frac{ \sqrt{K_{\max}}  }{h_{\min}} \norm{\bar{u}_h-\mathfrak{u}_h}_\Dom \leq C \frac{ \sqrt{K_{\max}}  }{h_{\min}} \norm{\bar{u}_h-v_h}_\Dom \leq C \frac{ \sqrt{K_{\max}}  }{h_{\min}}  \left( \norm{u-v_h}_\Dom + \norm{u-\bar{u}_h}_\Dom \right).
    \end{align*}
    For the jump term, the triangle inequality leads to 
    \begin{align*}
        \qquad\qquad&\sum_{F\in\calF_i\cup\calF_D} \frac{\sigma_F}{h_F}\norm{\jump{u-\mathfrak{u}_h}}_F^2 \leq 2\sum_{F\in\calF_i\cup\calF_D} \frac{\sigma_F}{h_F}\norm{\jump{u-\bar{u}_h}}_F^2 + 2\sum_{F\in\calF_i\cup\calF_D} \frac{\sigma_F}{h_F}\norm{\jump{\bar{u}_h-\mathfrak{u}_h}}_F^2.
    \end{align*}
    Then, with the help of the continuous trace inequality \eqref{eq_con_trace}, discrete inverse inequality \eqref{eq_dis_trace_inv}, $\sigma_F\leq \min\{\delta_F^+,\delta_F^-\}$ on $F\in\calF_i$ and $\sigma_F=\delta_F$ for $F\in\calF_D$, the triangle inequality and the fact that $\mathfrak{u}_h$ is the minimizer of \eqref{eq_scheme_step2_L2}, we derive
    \begin{align*}
        \qquad & \sum_{F\in\calF_i\cup\calF_D} \frac{\sigma_F}{h_F}\norm{\jump{\bar{u}_h-\mathfrak{u}_h}}_F^2 \leq C \sum_{T\in\calT_h} h_T^{-1}K_{\max,T} ( \norm{\nabla (\bar{u}_h-\mathfrak{u}_h)}_T + h_T^{-1}\norm{\bar{u}_h-\mathfrak{u}_h}_T )\norm{\bar{u}_h-\mathfrak{u}_h}_T\\
        &\leq C \sum_{T\in\calT_h} h_T^{-2}K_{\max,T} \norm{\bar{u}_h-\mathfrak{u}_h}_T^2 \leq C \frac{K_{\max}}{h_{\min}^2} \norm{\bar{u}_h-v_h}_{\Dom}^2 \leq C \frac{K_{\max}}{h_{\min}^2} (\norm{\bar{u}_h-u}_{\Dom}^2 + \norm{u-v_h}_{\Dom}^2).
    \end{align*}
    The above estimates imply
    \begin{equation*}
        \dgnorm{u-\mathfrak{u}_h}^2 \leq C \left( \dgnorm{u-\bar{u}_h}^2 + \left( \mu + \frac{K_{\max}}{h_{\min}^2} \right)\norm{u-v_h}_\Dom^2 + \frac{K_{\max}}{h_{\min}^2}\norm{u-\bar{u}_h}_\Dom^2 \right).
    \end{equation*}
    Owing to the regularity assumption, $u\in H^{p+1}(\calT_h)\subset H^2(\calT_h)$ for $p\geq 1$, the three terms on the right hand side are controlled by \eqref{eq_pre_H1_bound}, Lemma \ref{lemma_opt_vh_abm} and \eqref{eq_error_L2_pre}, respectively. This leads to the expected bound \eqref{eq_error_L2_cor_H1}.
\end{proof}

\subsection{A posteriori analysis}\label{sec_aposteriori_L2}

Let $\bar u_h$ and $\mathfrak u_h$ be given by \eqref{eq_scheme_step1} and \eqref{eq_scheme_step2_L2}. We define
\begin{equation}\label{eq_eta_L2}
    \eta_{L^2,T}^2\eqq\eta_{0,T}^2+\dgnorm{\mathfrak u_h-\bar u_h}_T^2,\qquad
    \eta_{L^2}^2\eqq\sum_{T\in\calT_h}\eta_{L^2,T}^2
    =\eta_0^2+\dgnorm{\mathfrak u_h-\bar u_h}^2.
\end{equation}

\begin{theorem}[A posteriori error estimate with $L^2$ correction]\label{thm_aposteriori_L2}
    Under the assumptions of Lemma \ref{lemma_reliability_sta}, $\mathfrak u_h$ satisfies
    \begin{equation}\label{eq_reliability_bound_L2}
        \dgnorm{u-\mathfrak u_h}^2\leq C\eta_{L^2}^2,
    \end{equation}
    and
    \begin{equation}\label{eq_efficiency_bound_L2}
        \eta_{L^2}^2\leq C(\mu,K)\big(\dgnorm{u-\mathfrak u_h}^2+\calR_{osc}^2\big).
    \end{equation}
    Moreover,
    \begin{equation}\label{eq_correction_efficiency_L2}
        \dgnorm{\mathfrak u_h-\bar u_h}^2
        \leq C(\mu,K)\big(\dgnorm{u-\mathfrak u_h}^2+\calR_{osc}^2\big),
    \end{equation}
    and, for every $T\in\calT_h$,
    \begin{equation}\label{eq_local_efficiency_correction_L2}
        \eta_{L^2,T}^2\leq C(\mu,K)\sum_{T'\in\calT_T}
        \big(\dgnorm{u-\mathfrak u_h}_{T'}^2+\dgnorm{\mathfrak u_h-\bar u_h}_{T'}^2+\calR_{osc,T'}^2\big).
    \end{equation}
\end{theorem}

\begin{proof}
    Since $\mathfrak u_h\in V_h^p$, Corollary \ref{cor_arbitrary_correction} with $v_h=\mathfrak u_h$ proves \eqref{eq_reliability_bound_L2} and \eqref{eq_efficiency_bound_L2}. The estimate \eqref{eq_correction_efficiency_L2} follows from \eqref{eq_prediction_comparison} with the same choice. Finally, \eqref{eq_prediction_local_efficiency} and
    \begin{equation*}
        \dgnorm{u-\bar u_h}_{T'}^2
        \leq2\dgnorm{u-\mathfrak u_h}_{T'}^2
        +2\dgnorm{\mathfrak u_h-\bar u_h}_{T'}^2
    \end{equation*}
    give \eqref{eq_local_efficiency_correction_L2}.
\end{proof}

\begin{remark}[extension to other prediction-correction schemes]\label{rmk_other_correction}
    The above proof uses only that the corrected approximation belongs to $V_h^p$. Consequently, the same a posteriori analysis applies to any correction of the predictor $\bar{u}_h$, with the estimator defined in Corollary \ref{cor_arbitrary_correction}.
\end{remark}

\section{Numerical experiments}\label{sec_num_test}
In this section, we numerically test and compare the proposed schemes \eqref{eq_scheme_step1}-\eqref{eq_scheme_step2} and \eqref{eq_scheme_step1}-\eqref{eq_scheme_step2_L2}. The schemes are implemented with the help of \texttt{deal.ii} library \cite{dealii_2025}.

We consider three numerical test cases in this section. The first one is an accuracy test with a piecewise smooth exact solution, and mixed Dirichlet and Neumann boundary conditions are imposed; the second and the third one are tests with anisotropic diffusion, and the closed form of the exact solutions is unknown. In addition, the second test case has a mildly anisotropic but heterogeneous diffusion, whereas the third test case has strongly anisotropic diffusion.

\paragraph*{Test case 1.} Here we consider the PDE  $u-\DIV(K\nabla u) = f$ in the domain $\Dom=(0,1)^2$, with the following exact solution and parameters:
\begin{equation}\label{eq_test1}
    K=\begin{cases}
        \text{diag}(1,\epsilon) &\text{ if } x<0.5\\
        \text{diag}(\epsilon^{-1},1) &\text{ if } x>0.5\\
    \end{cases},\quad 
    u(x,y)=\begin{cases}
        \frac{2(1+\cos(\pi y))}{\epsilon(1+\epsilon^{-1})} x^2 &\text{ if } 0\leq x \leq 0.5\\
        \frac{(1+\cos(\pi y)) (1+\epsilon^{-1}-4(1-x)^2) }{2(1+\epsilon^{-1})}  &\text{ if } 0.5\leq x \leq 1
    \end{cases},
\end{equation}
and $\epsilon=10^{-2}$, $f$ is constructed from above functions, Neumann boundary conditions are imposed on left and bottom boundary, and Dirichlet boundary conditions are imposed on the rest of the boundary. We notice that such $u$ is piecewise smooth, hence dG methods are expected to have optimal convergence on meshes perfectly matching the interface at $\{x=0.5\}$. We have $a=0$ and $b=1$ in this test case.

\paragraph*{Test case 2.} We consider the PDE  $u-\DIV(K\nabla u) = f$ in the domain $\Dom=(0,1)^2$, with the following parameters: {\color{red} Define $R_{\theta}$ below!}
\begin{subequations}\label{eq_test2}
\begin{align}
    &K(x,y)=R_{\theta}\begin{pmatrix}
        \lambda_1(x,y) & 0 \\
        0 & \lambda_2(x,y) 
    \end{pmatrix}R_{\theta}^{-1}, \quad 
    f(x,y)=\begin{cases}
        1.0 &\text{ if } (x,y)\in (0,0.5)\times (0,1)\\
        0.5 &\text{ if } (x,y)\in (0.5,1)\times (0.5,1)\\
        0 &\text{ if } (x,y)\in  (0.5,1)\times (0,0.5),
    \end{cases},\\
    &\lambda_1(x,y)=\begin{cases}
        5\times 10^{-3} &\text{ if } 0<x<0.5\\
        1  &\text{ if } 0.5<x<1\\
    \end{cases}, \quad 
    \lambda_2(x,y)=\begin{cases}
        1 &\text{ if } 0<x<0.5\\
        5\times 10^{-3}  &\text{ if } 0.5<x<1
    \end{cases},
\end{align}
\end{subequations}
and $\theta=40^\circ$, $K\nabla u \cdot \bn = 0$ on the boundary $\partial \Dom$. Interfaces are expected near the discontinuities of data in this test case. We have $a=0$ and $b=1$ in this test case.

\paragraph*{Test case 3.} We consider the following PDE inspired by \cite[equation (42)]{Kirby_Shapero_2024} :  $u-\DIV(K\nabla u) = f$ in the domain $\Dom=(0,1)^2$, with parameters
\begin{equation}\label{eq_test3}
    K(x,y)=\begin{bmatrix}
        x^2+\epsilon y^2 & -(1-\epsilon)xy \\
        -(1-\epsilon)xy & y^2+\epsilon x^2 
        \end{bmatrix}, \qquad f(x,y)=\begin{cases}
            1, \text{ if } |(x,y)-(\frac{1}{2},\frac{1}{2})|\leq \frac{1}{8}\\
            0, \text{ otherwise}
        \end{cases},
\end{equation}
$\epsilon=10^{-4}$, and the homogeneous Dirichlet boundary condition is imposed for $u$. This test has a variable, degenerate diffusion tensor and lies outside the coefficient assumptions of the analysis. It is included to examine the numerical behavior of the proposed corrections in a more challenging setting.

There are two numerical difficulties in this test case: small eigenvalue of $O(\epsilon)$ of diffusion tensor $K$ and $K=0$ at $(0,0)$. The maximum principle implies that $u$ is positive in this test case; for the upper bound, we set a rough bound $b=1$, since a more accurate upper bound is not known a priori. In addition, $K$ is not piecewise constant, in practical simulations, $K$ and its related quantities are evaluated at quadrature points which do not contain $(0,0)$.

We examine the performance of the proposed schemes from four perspectives. First, we test their accuracy using Test case 1. Second, we compare the corrected solutions \eqref{eq_scheme_step2} and \eqref{eq_scheme_step2_L2} with the prediction step \eqref{eq_scheme_step1} to illustrate the effect of enforcing the constraints. Third, we investigate adaptive refinement for Test cases 2 and 3. Finally, we compare our estimators with those derived using variational inequality techniques.

Moreover, we list a few numerical conventions used in this section:

\paragraph*{DG stabilization parameter.} Although the $p$-dependence is not investigated in this work, we use the standard $p^2$-scaled penalty choice $\gamma_0=10p(p+1)$.

\paragraph*{Meshes.} All meshes are produced by default setting of \texttt{deal.ii} library, hence for the uniform refinement, we use the structured tensor-product meshes; for the adaptive refinement, meshes are not conforming, in the sense that each edge has at most one hanging node. We note that the nonconforming rectangle meshes is not covered by the current theory, but allows a simplicial shape-regular sub-triangulation. 

\paragraph*{PDAS solver.} In all our numerical experiments, the PDAS converges within 5 iterations, hence we decide not to report its iteration numbers in this section, for the sake of brevity. 

\paragraph*{Structure-preserving set $X_h$.} We only consider the case $X_h=X_h^{0,1,m}$, since this set is more challenging than the bound-preserving case, and we observe essentially the same convergence behavior for both $X_h$. For Test case 1, we set the total mass by $m=\int_\Dom u$; for Test case 2, we naturally set $m=\mu^{-1}(\int_\Dom f + \int_{\partial \Dom} g_N)$; for Test case 3, we set $m=\int_\Dom \bar{u}_h$, hence the correction preserves the prediction step’s total mass, which is not mathematically rigorous, but practically works. 

\paragraph*{Bound-violation.} In all numerical tests, the prediction step \eqref{eq_scheme_step1} produces nodal values outside the prescribed bounds, while both correction steps restore nodal bound-preservation. For brevity, we do not report these violations and their correction separately.

\paragraph*{Effectivity index.}
Effectivity index is commonly used to measure the sharpness of the proposed estimators, we define it by
\begin{equation}
    \text{effectivity index}(v_h) = \frac{\eta}{\dgnorm{u-v_h}} \quad v_h \in \{ \bar{u}_h, u_h, \mathfrak{u}_h\}.
\end{equation}

\paragraph*{Adaptive scheme.} We consider here an adaptive algorithm
with the usual structure
\begin{equation}\label{eq_adaptive_scheme}
    \textbf{SOLVE}
    \qquad \Longrightarrow \qquad
    \textbf{ESTIMATE}
    \qquad \Longrightarrow \qquad
    \textbf{MARK}
    \qquad \Longrightarrow \qquad
    \textbf{REFINE}.
\end{equation}
This structure is iteratively applied until either a given tolerance for $\eta$ is achieved or a given maximum iteration number $N_{\rm{iter}}$ is attained.
The \textbf{ESTIMATE} step is driven by using the local error estimator $\eta|_T$ on each mesh cell $T\in \calT_h$. 
As for the \textbf{MARK} step, we use D\"orfler’s marking
with a given threshold $\theta_{\rm{refine}}$ in~$(0,1]$ \cite{Dorfler_1996}.
The \textbf{REFINE} step is realized by the default refine strategy of the \texttt{deal.ii} library.

\paragraph*{Data abbreviation.} Data associated to numerical solutions $\bar{u}_h$ from \eqref{eq_scheme_step1}, $u_h$ from \eqref{eq_scheme_step2} and $\mathfrak{u}_h$ from \eqref{eq_scheme_step2_L2} are named with prefix `no corr.', `dG corr.'  and `L2 corr.', respectively.

\subsection{Accuracy test on a piecewise smooth test case}
\begin{figure}[!htbp]
    \centering
    \includegraphics[width=0.3\linewidth]{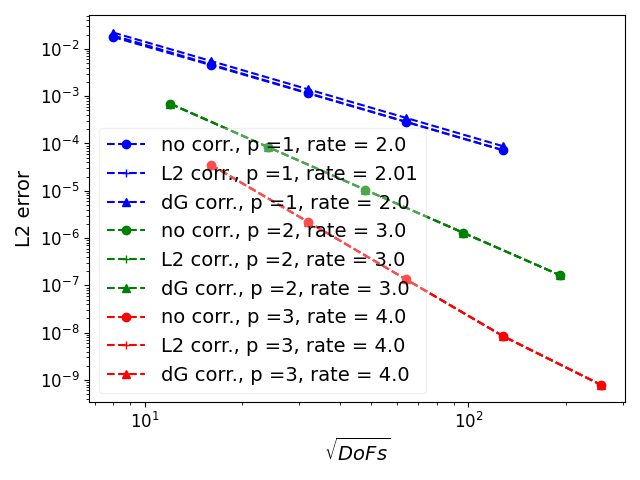}
    \includegraphics[width=0.3\linewidth]{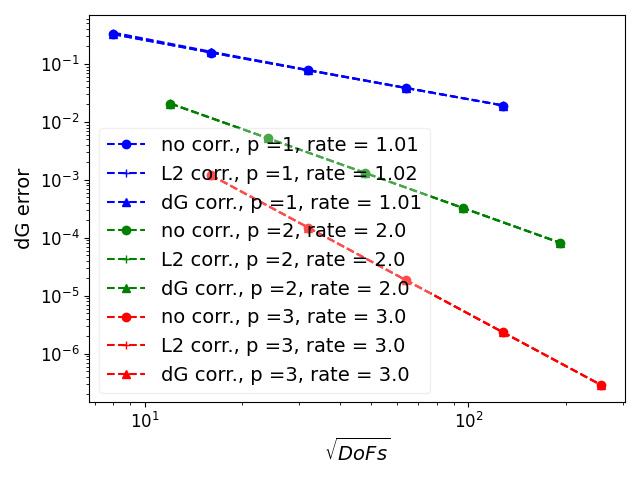}
    \includegraphics[width=0.3\linewidth]{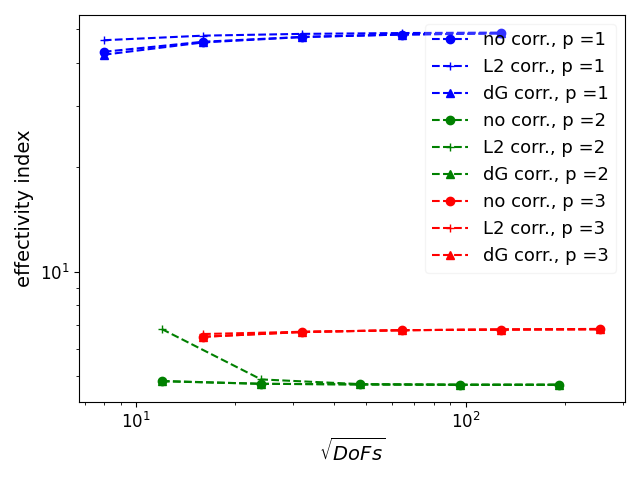}
    \caption{Accuracy test on \eqref{eq_test1}: left: $L^2$ error; middle: dG energy error; right: effectivity index.}
    \label{fig:accuracy}
\end{figure}

In Figure \ref{fig:accuracy}, we illustrate the numerical results obtained with $p\in\{1,2,3\}$ and uniform mesh with $h\in\{2^{-2},2^{-3},2^{-4},2^{-5},2^{-6}\}$, for the test case \eqref{eq_test1}. The $L^2$-, dG energy errors and effectivity index are reported for the prediction step \eqref{eq_scheme_step1}, scheme \eqref{eq_scheme_step1}-\eqref{eq_scheme_step2} and \eqref{eq_scheme_step1}-\eqref{eq_scheme_step2_L2}.

From Figure \ref{fig:accuracy}, we have the following observations: 1. the optimal convergence rates are observed for the $L^2$ norm, dG energy norm and error estimator $\eta$; 2. the effectivity index is uniformly bounded for any fixed $p$; 3. the error curves with and without correction almost overlap, indicating that the correction steps enforce the constraints without degrading convergence.

\subsection{Behavior of two correction steps.}
\begin{figure}[!htbp]
    \centering

    \begin{subfigure}{0.26\textwidth}
    \centering
    \includegraphics[width=\textwidth]{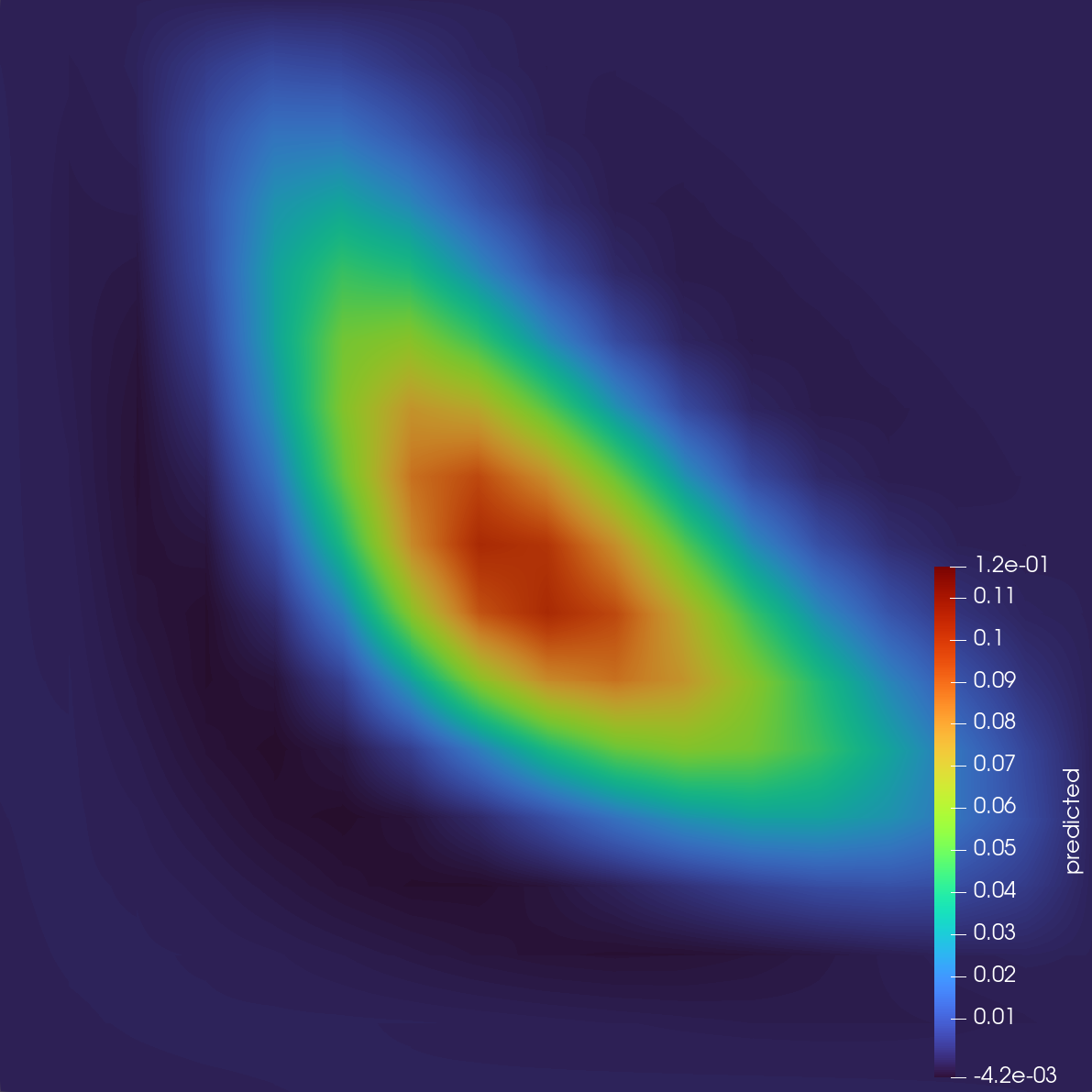}
    \caption{Pred. sol., coarse.}
    \end{subfigure}
    \begin{subfigure}{0.26\textwidth}
    \centering
    \includegraphics[width=\textwidth]{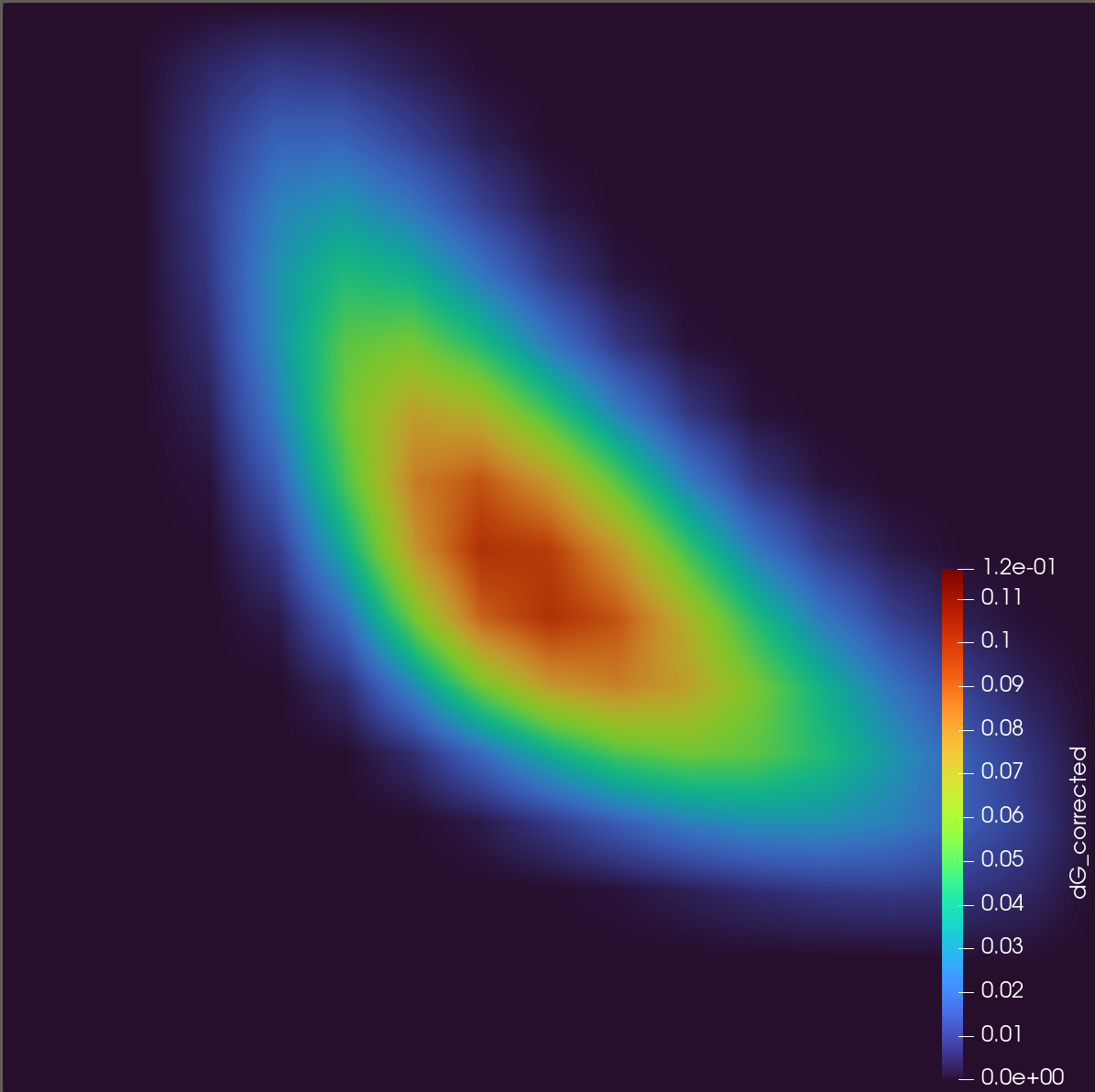}
    \caption{DG corr. sol., coarse.}
    \end{subfigure}
    \begin{subfigure}{0.26\textwidth}
    \centering
    \includegraphics[width=\textwidth]{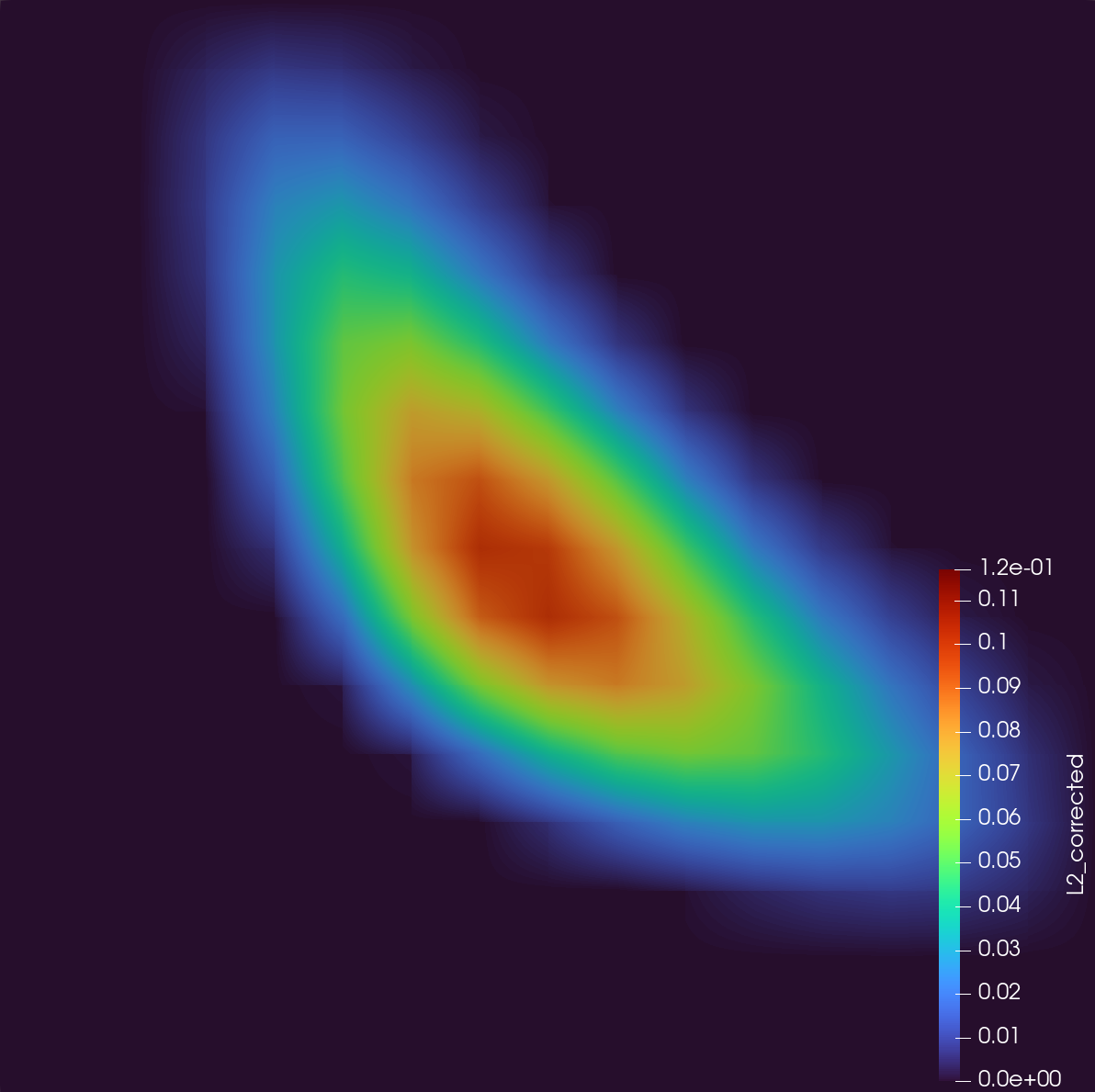}
    \caption{$L^2$ corr. sol., coarse.}
    \end{subfigure}

    \begin{subfigure}{0.26\textwidth}
    \centering
    \includegraphics[width=\textwidth]{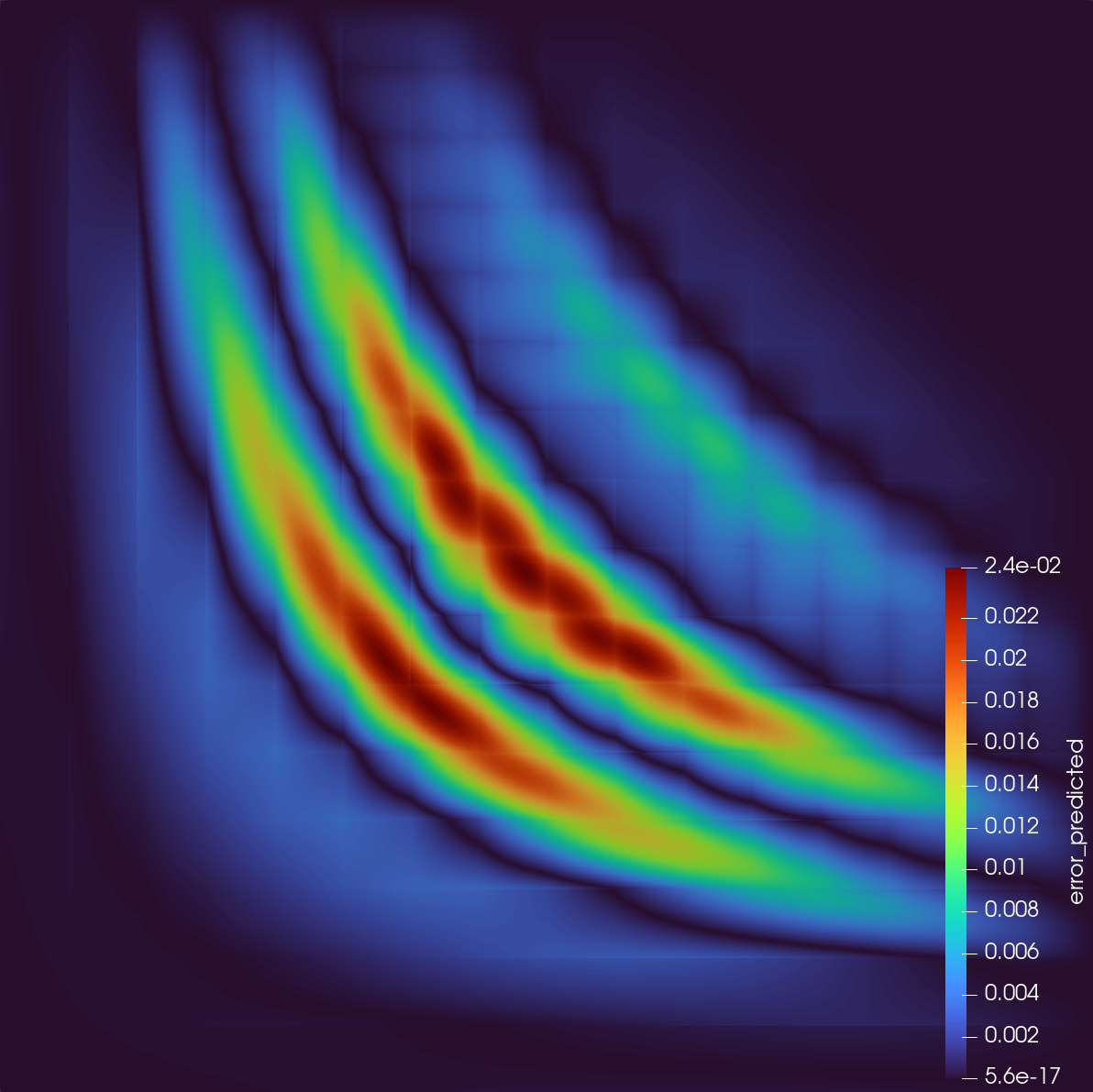}
    \caption{Error, pred., coarse.}
    \end{subfigure}
    \begin{subfigure}{0.26\textwidth}
    \centering
    \includegraphics[width=\textwidth]{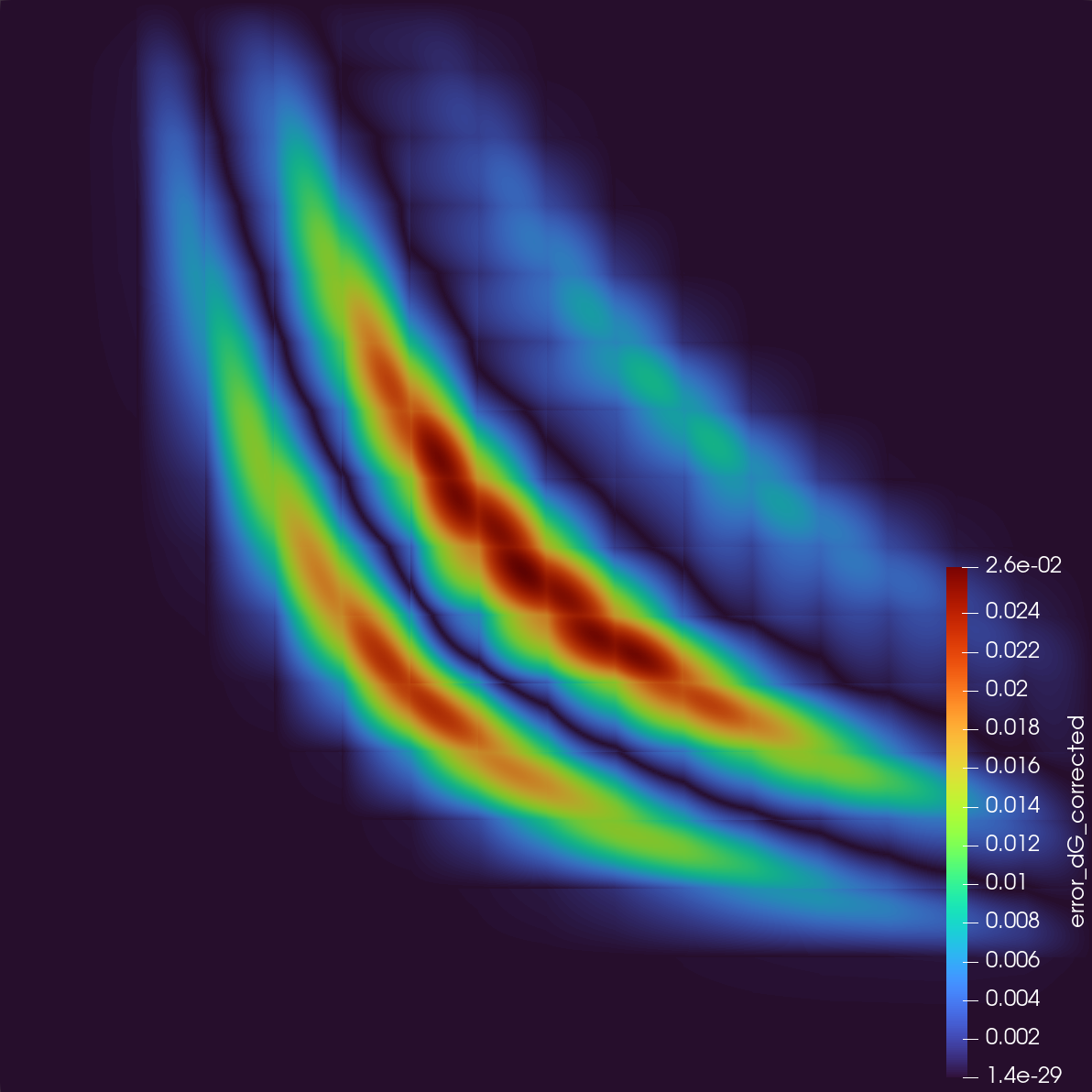}
    \caption{Error, dG corr., coarse.}
    \end{subfigure}
    \begin{subfigure}{0.26\textwidth}
    \centering
    \includegraphics[width=\textwidth]{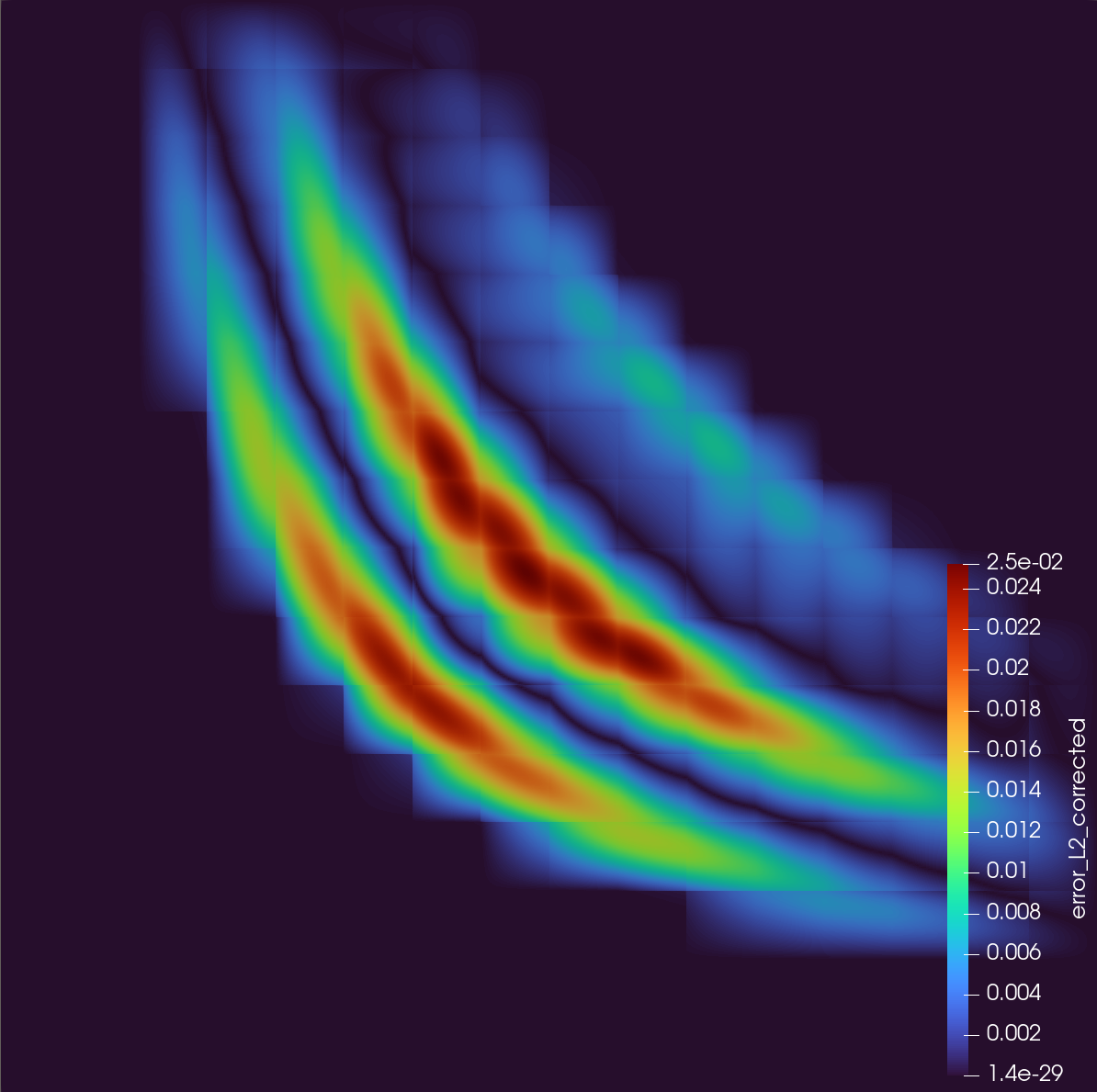}
    \caption{Error, $L^2$ corr., coarse.}
    \end{subfigure}

    \begin{subfigure}{0.26\textwidth}
    \centering
    \includegraphics[width=\textwidth]{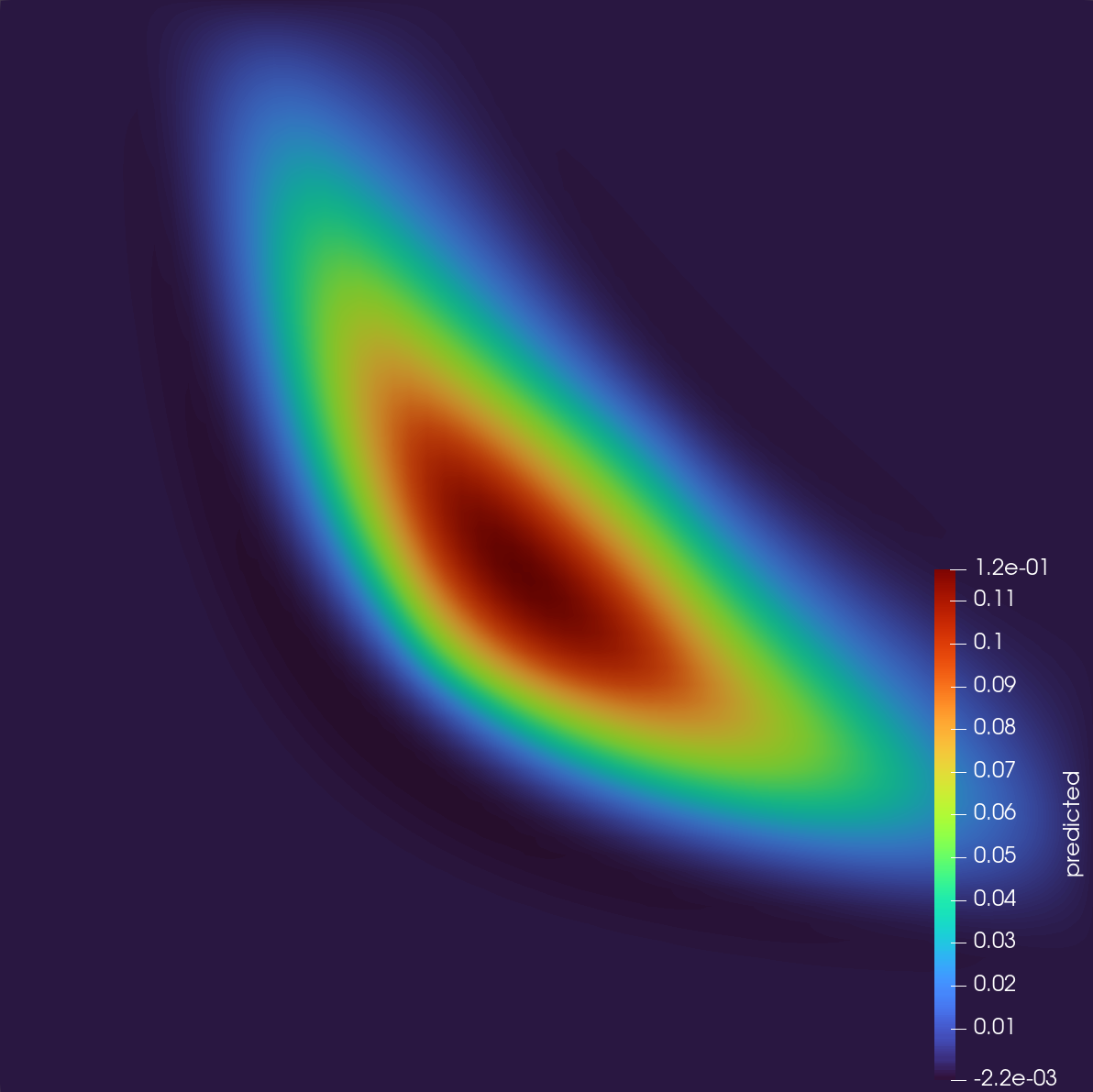}
    \caption{Pred. sol., fine.}
    \end{subfigure}
    \begin{subfigure}{0.26\textwidth}
    \centering
    \includegraphics[width=\textwidth]{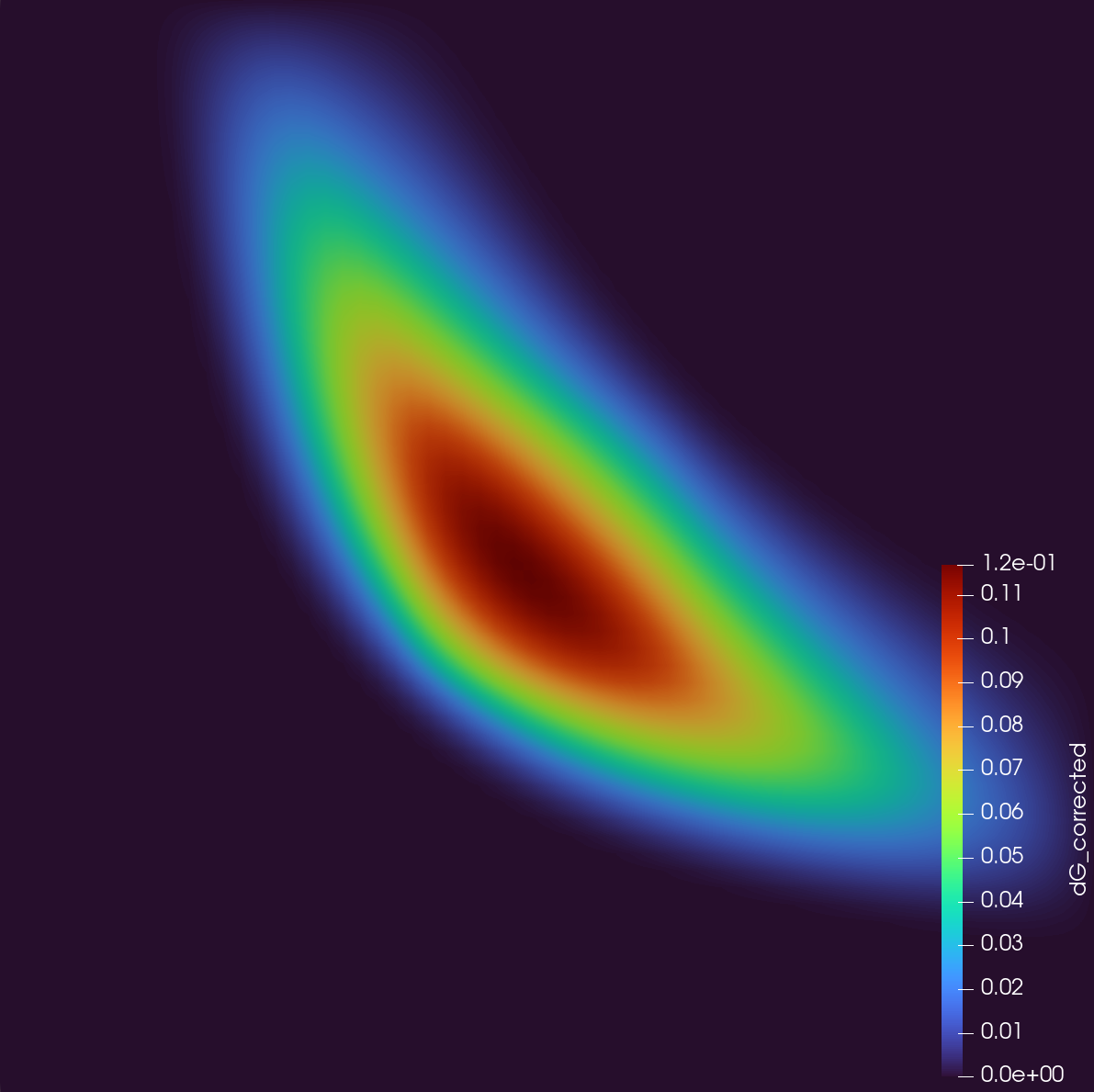}
    \caption{DG corr. sol., fine.}
    \end{subfigure}
    \begin{subfigure}{0.26\textwidth}
    \centering
    \includegraphics[width=\textwidth]{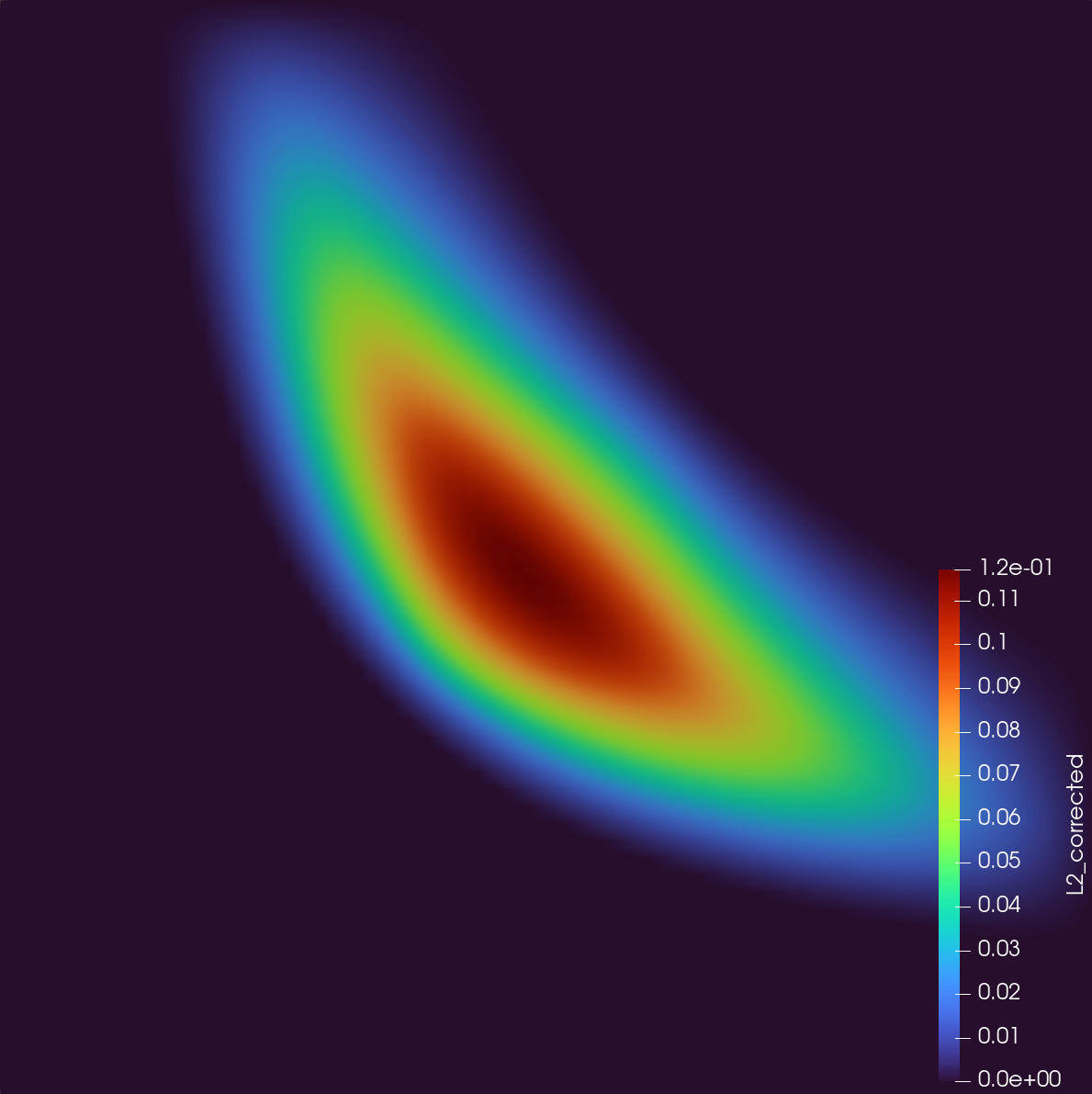}
    \caption{$L^2$ corr. sol., fine.}
    \end{subfigure}

    \begin{subfigure}{0.26\textwidth}
    \centering
    \includegraphics[width=\textwidth]{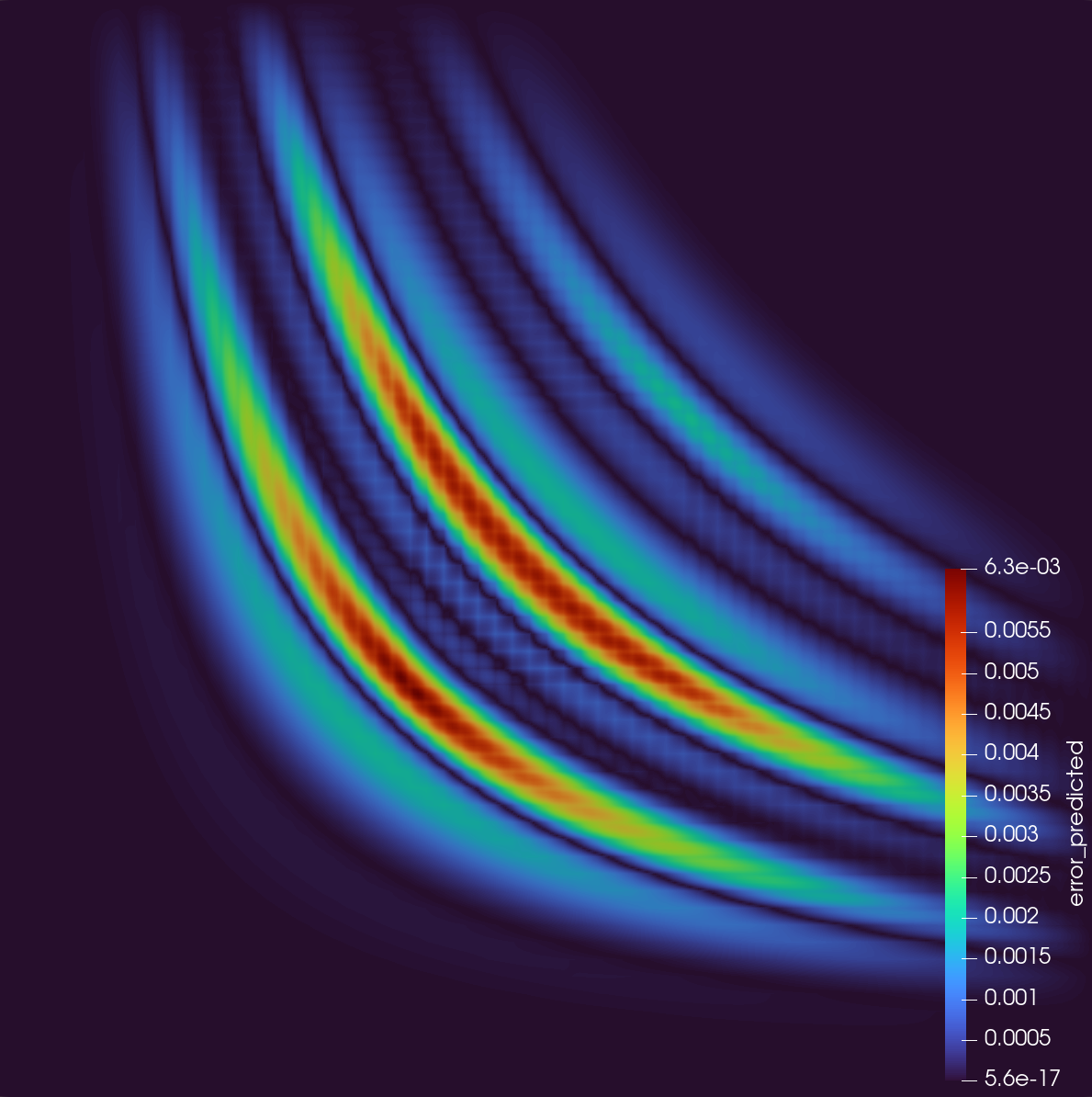}
    \caption{Error, pred., fine.}
    \end{subfigure}
    \begin{subfigure}{0.26\textwidth}
    \centering
    \includegraphics[width=\textwidth]{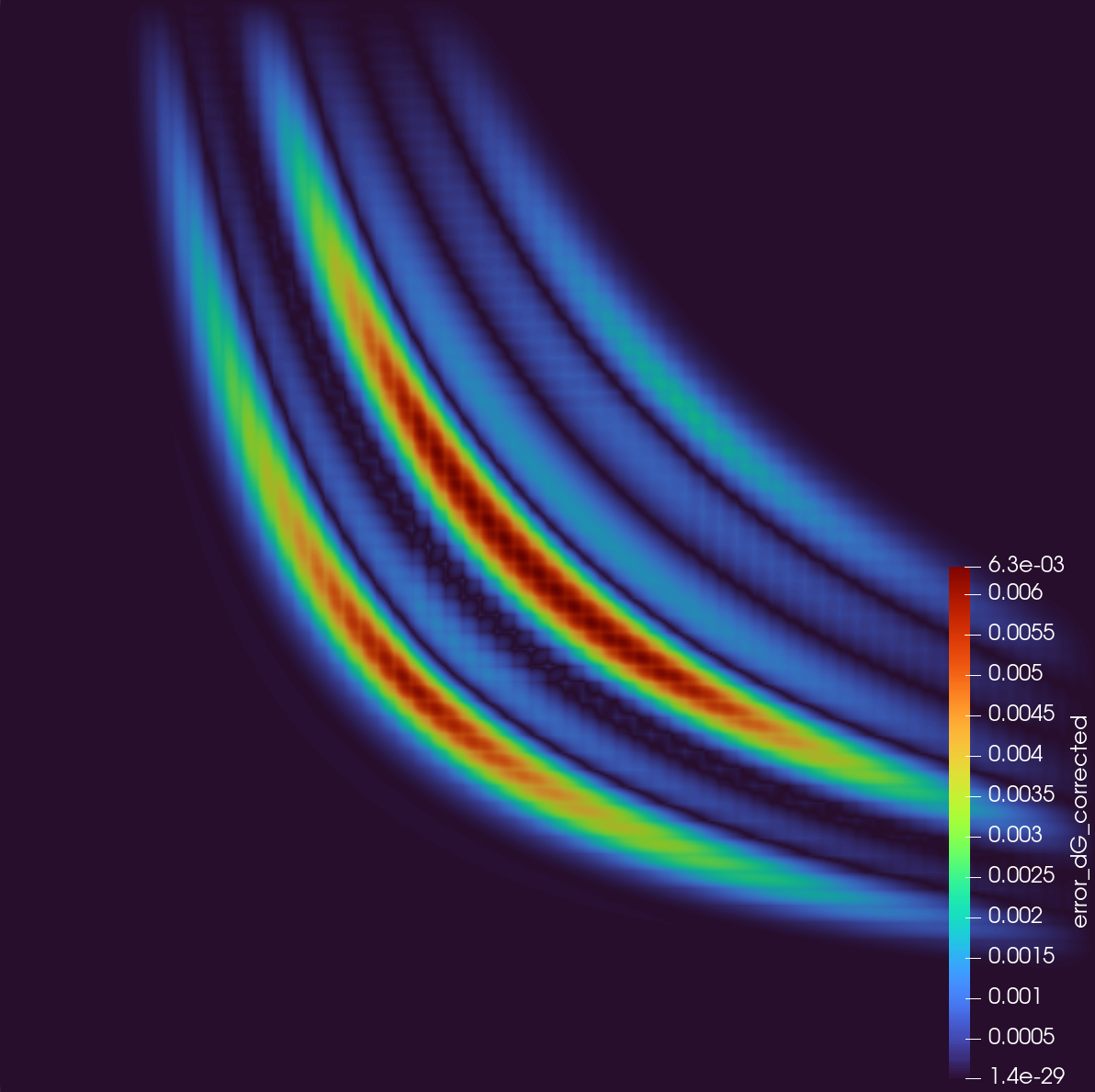}
    \caption{Error, dG corr., fine.}
    \end{subfigure}
    \begin{subfigure}{0.26\textwidth}
    \centering
    \includegraphics[width=\textwidth]{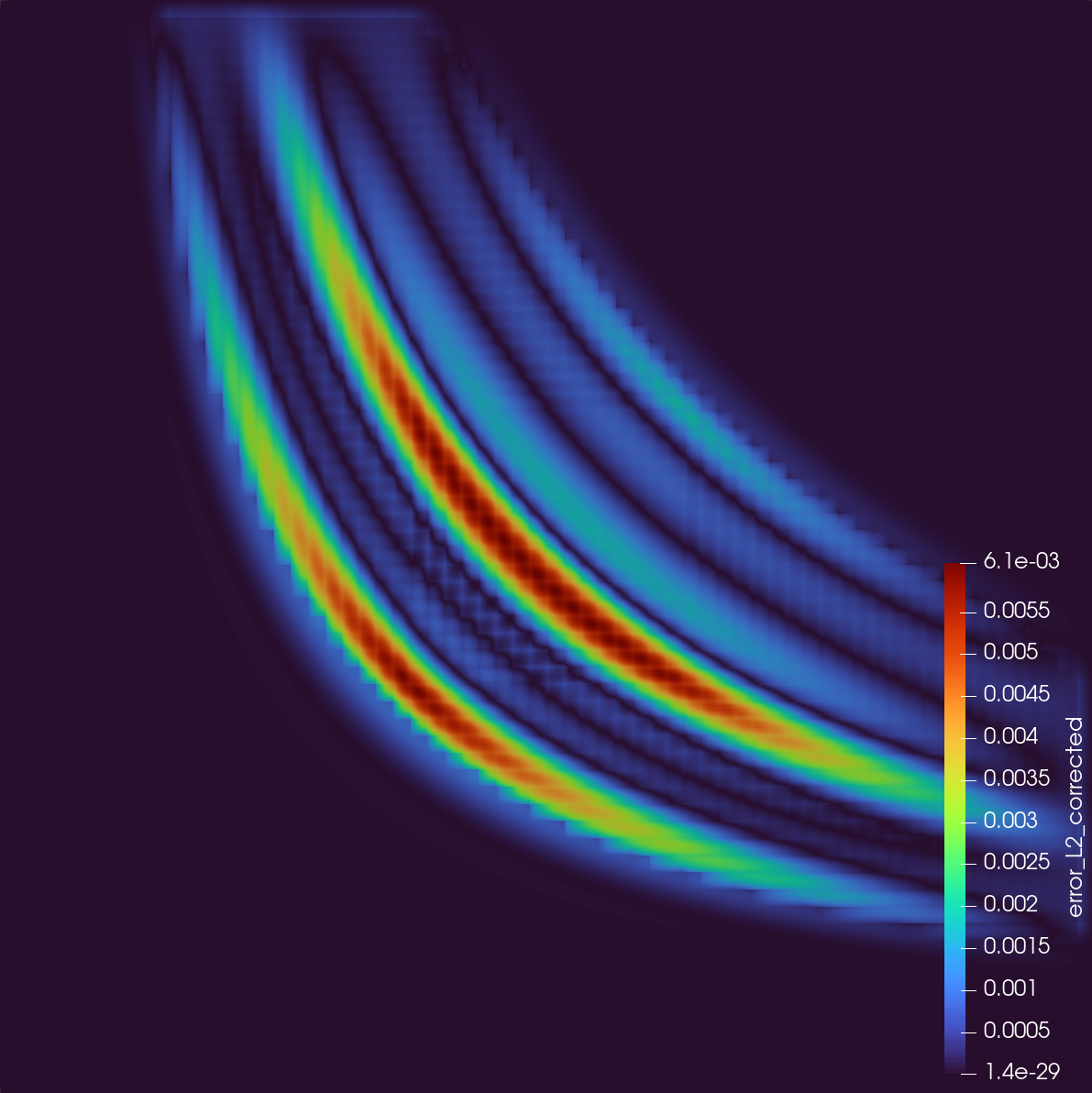}
    \caption{Error, $L^2$ corr., fine.}
    \end{subfigure}
    
    \caption{ Numerical solution for test case 3 in \eqref{eq_test3}: $p=1$, reference, coarse and fine meshes are generated with $N=256,16,64$ intervals in each direction.}
    \label{fig:error_Kirby}
\end{figure}

We illustrate the effect of correction steps by considering Test case 3, on coarse and fine meshes, respectively.
To fairly illustrate and compare the effect of correction steps, we visualize the numerical solution with/without correction, and their pointwise error plots in Figure \ref{fig:error_Kirby}, 
with a numerical reference solution $u_{ref}$ produced with $p=1$ and $N=256$ intervals in each direction. The numerical solutions used for comparison are produced with $p=1$ and $N=16,64$.

From Figure \ref{fig:error_Kirby}, we observe that: 1. the negative parts of the prediction steps are well corrected for both correction steps, even for the coarse mesh; 2. the dG energy based correction leads to a smoother solution and a smoother error distribution near the transition region $\{\bar{u}_h \approx 0\}$, for the coarse mesh, and such behavior is less evident for the fine mesh. These observations suggest using the dG energy (resp. $L^2$) correction for underresolved (resp. resolved) diffusion.

\subsection{Performance of adaptive scheme}\label{sec_adaptive_scheme}
\begin{figure}[!htbp]
    \centering

    \includegraphics[width=0.45\linewidth]{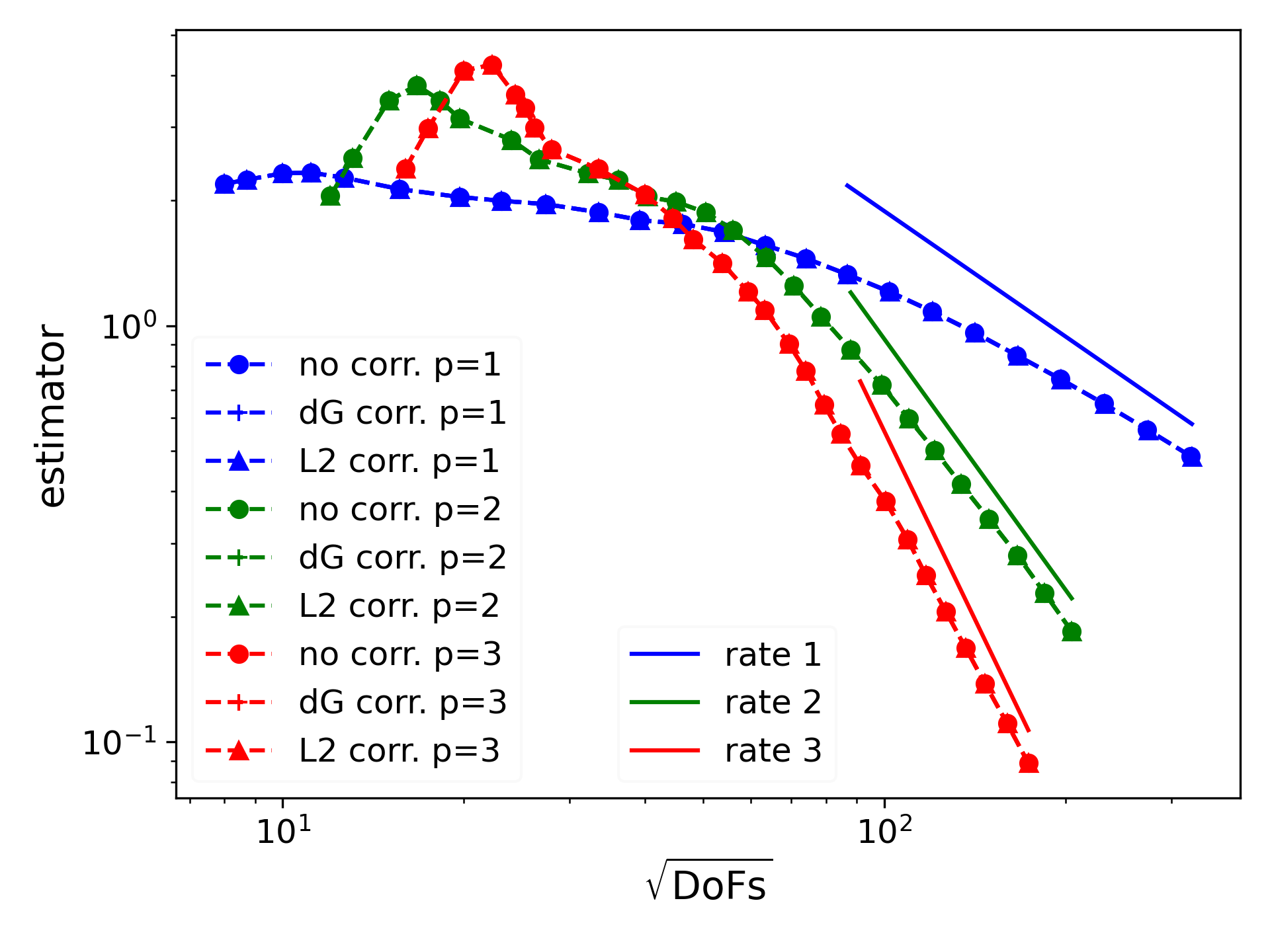}
    \includegraphics[width=0.45\linewidth]{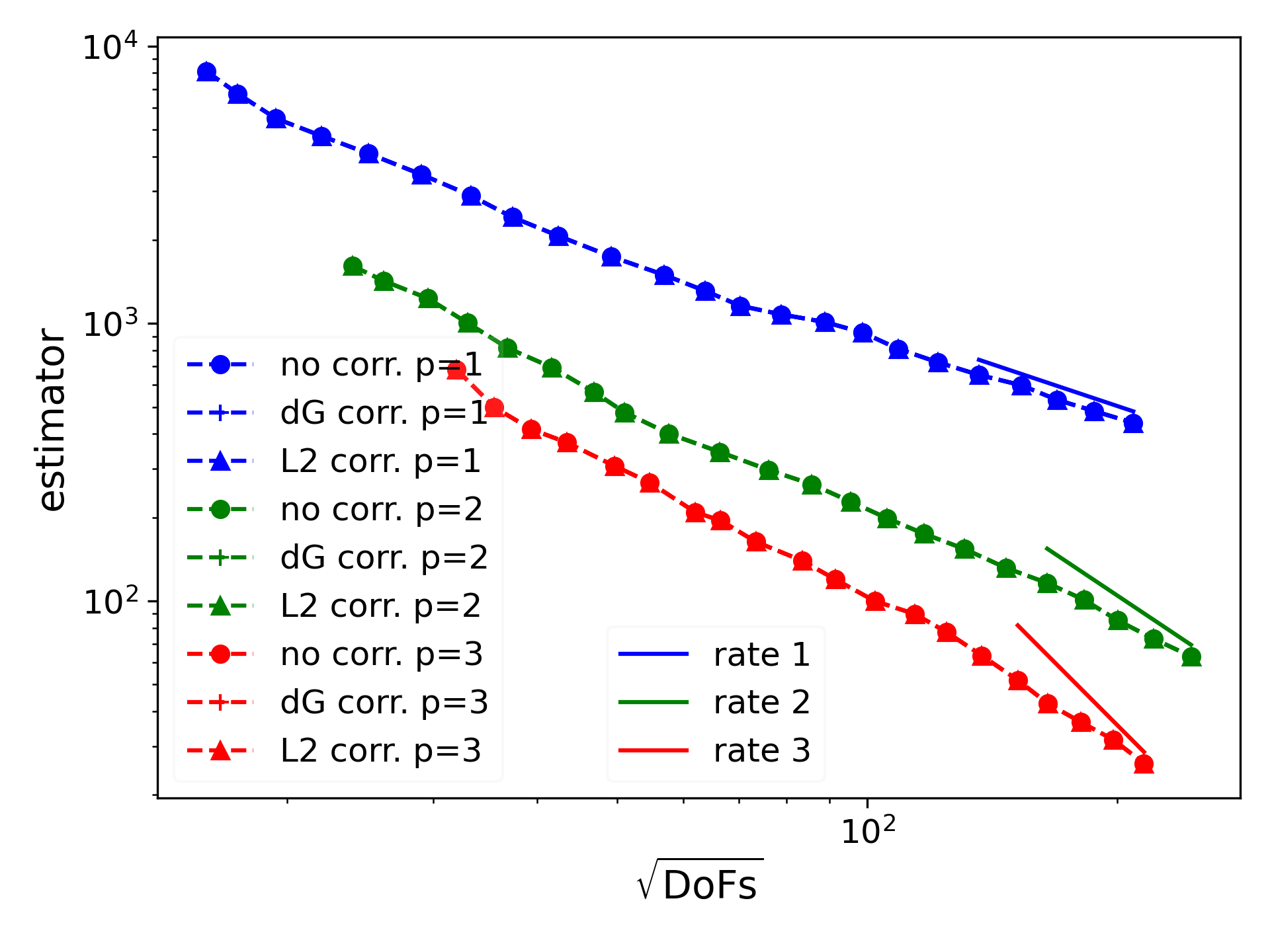}
    
    \caption{Error estimators with adaptive scheme \eqref{eq_adaptive_scheme}: left: Test case 2; right: Test case 3.}
    \label{fig:adaptive_rate}
\end{figure}

In this part, we test the adaptive scheme \eqref{eq_adaptive_scheme} for prediction step and two correction steps, for Test cases 2 and 3.
In Figure \ref{fig:adaptive_rate}, we set $\theta_{refine}=0.3$ and uniform initial mesh with $h=2^{-2}$, and present the estimators for $p\in\{1,2,3\}$ for schemes \eqref{eq_scheme_step1}, \eqref{eq_scheme_step1}-\eqref{eq_scheme_step2} and \eqref{eq_scheme_step1}-\eqref{eq_scheme_step2_L2}. 

We observe from Figure \ref{fig:adaptive_rate} that: for both the mildly and strongly anisotropic test cases, it shows that the estimator curves for both correction steps overlap those for the prediction scheme \eqref{eq_scheme_step1}. These results support using the less expensive $L^2$ correction. The adaptive scheme also identifies singular regions and refines the mesh locally.

\begin{figure}[!htbp]
    \centering
    \includegraphics[width=0.45\linewidth]{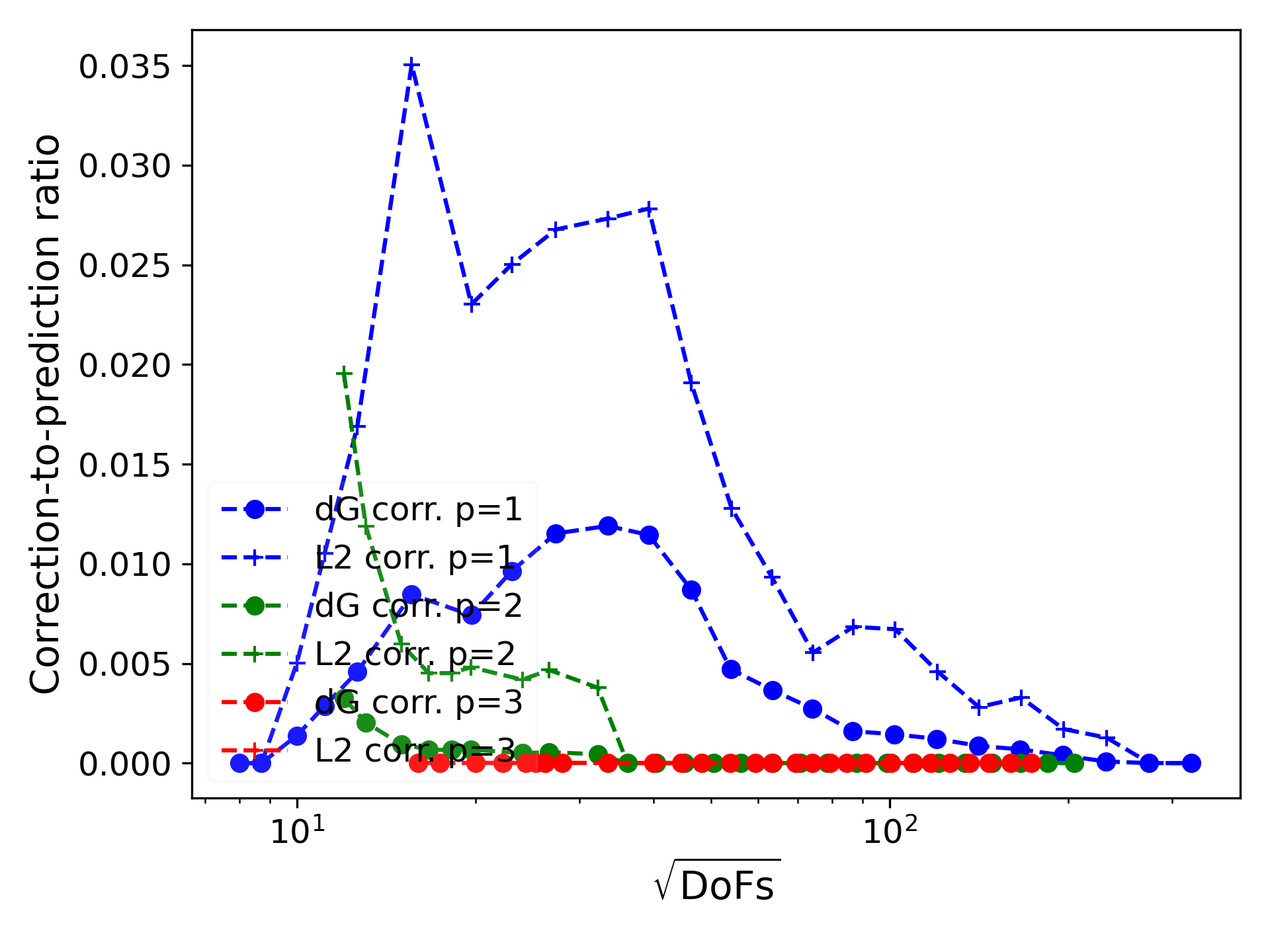}
    \includegraphics[width=0.45\linewidth]{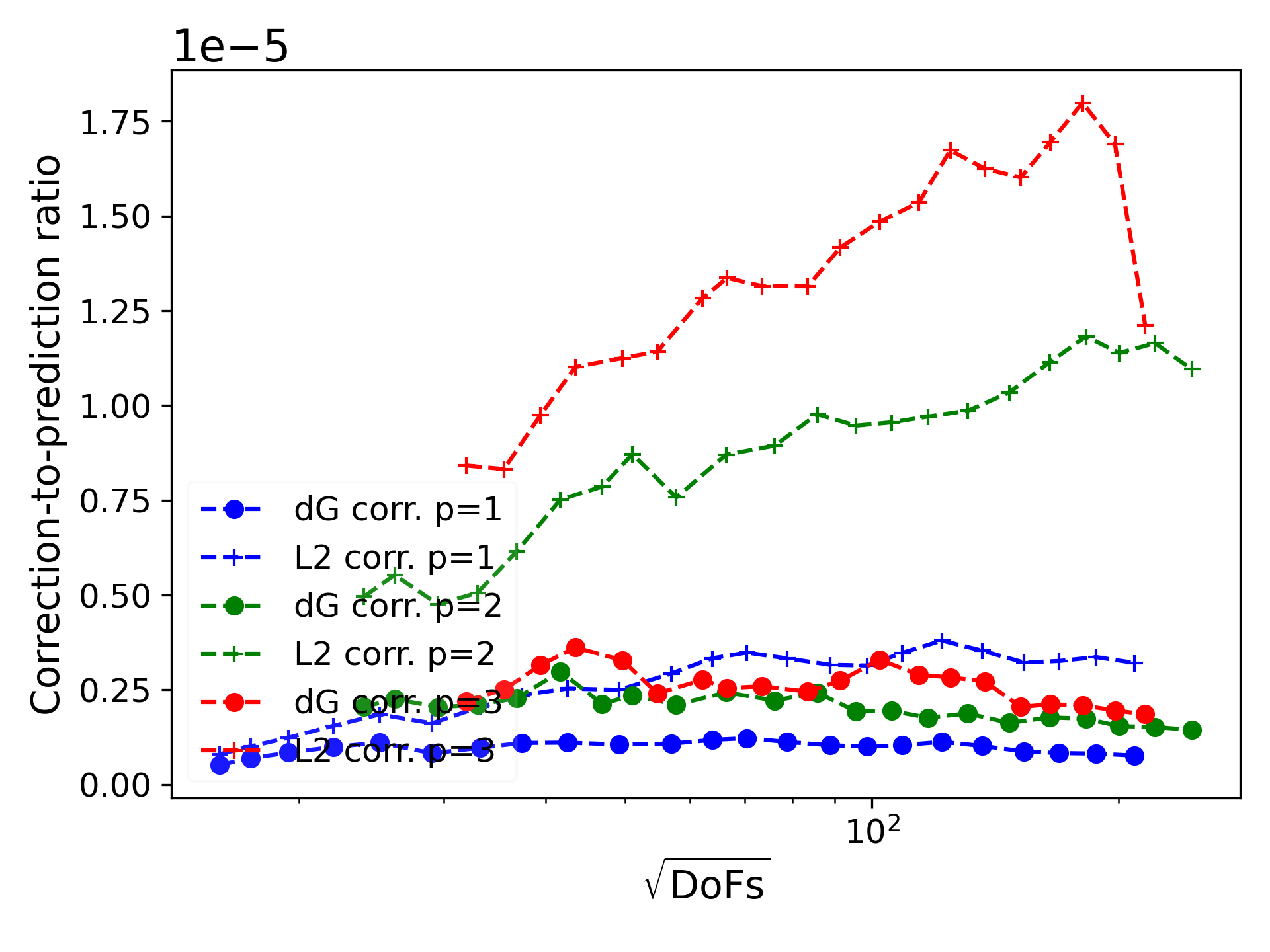}
    \caption{Prediction-correction estimators ratio with adaptive scheme \eqref{eq_adaptive_scheme}: left: test case 2 in \eqref{eq_test2}; right: Test case 3.}
    \label{fig:adaptive_ratio}
\end{figure}

Next, we explore the ratio between prediction and correction estimators, \ie $\dgnorm{u_h-\bar{u}_h}/ \eta_0$. In Figure \ref{fig:adaptive_ratio}, we illustrate this ratio for Test cases 2 and 3, with the same parameters used in Figure \ref{fig:adaptive_rate}.

From Figure \ref{fig:adaptive_ratio}, we observe that the $L^2$ correction leads to relatively larger correction estimator than the dG energy correction, but both of them are negligible compared to the prediction estimators. These results also imply that $\eta_0$ plays the main role in the marking step of the adaptive schemes, which explains the overlapped curves in Figure \ref{fig:adaptive_rate}. This observation suggests the usage of the $L^2$ correction in practice, since it is cheaper.

\subsection{Comparison with variational inequality argument}

The main goal of this section is to numerically compare our estimators and the estimators used in the variational inequality literature \cite{Braess_2005,Veeser_2001}.

We notice that with the help of $\sigma_h$ from \eqref{eq_sigma_h}, the numerical scheme with dG energy correction \eqref{eq_scheme_step1}-\eqref{eq_scheme_step2} can be reformulated as
\begin{equation}
    a_h(u_h,v_h) =l_h(v_h) + (\sigma_h,v_h) \quad \forall v_h \in V_h^p.
\end{equation}
Therefore, applying the a posteriori techniques used in the variational inequality literature \cite{Braess_2005,Veeser_2001}, it can be shown that
\begin{align}\label{eq_estimator_LM_dG_energy}
    \dgnorm{u-u_h}^2 &\leq C \Big[ \sum_{T\in\calT_h} \Big( \alpha_{1,T}^2 \norm{(I-\Pi_0) (f + \DIV (K\nabla_h u_h) -\mu u_h +\sigma_h)  }_T^2 \\
    &\quad + \eta_{J,T}(u_h)^2 + \eta_{N,T}(u_h)^2 + \chi_K ( \eta_{i,T}(u_h)^2 + \eta_{D,T}(u_h)^2 ) \Big) + \Big( \sup_{v\in V_0, v\neq 0}\frac{(\sigma_h,v)_{\Dom}}{\norm{v}_{E,\Dom}} \Big)^2 \Big].\nonumber
\end{align}
A corresponding efficiency bound can also be derived. The main differences between this variational inequality based estimator and our estimator are the appearance of the dual norm of Lagrange multiplier and the modified volumetric residual term.

Recall that, for the $L^2$ correction based scheme \eqref{eq_scheme_step1}-\eqref{eq_scheme_step2_L2}, we can reformulate it as \eqref{eq_scheme_L2_LM}, hence, following the variational inequality argument, we can find
\begin{align}\label{eq_estimator_LM_L2}
    \dgnorm{u-\mathfrak{u}_h}^2 &\leq C \Big[ \sum_{T\in\calT_h} \Big( \alpha_{1,T}^2 \norm{(I-\Pi_0) (f + \DIV (K\nabla_h \mathfrak{u}_h) -\mu \mathfrak{u}_h +\calL_h\sigma_{h,0})  }_T^2 \\
    &\quad + \eta_{J,T}(\mathfrak{u}_h)^2 + \eta_{N,T}(\mathfrak{u}_h)^2 + \chi_K( \eta_{i,T}(\mathfrak{u}_h)^2 + \eta_{D,T}(\mathfrak{u}_h)^2 ) \Big) + \Big( \sup_{v\in V_0, v\neq 0}\frac{(\calL_h\sigma_{h,0},v)_{\Dom}}{\norm{v}_{E,\Dom}} \Big)^2 \Big].\nonumber
\end{align}
A corresponding efficiency bound can also be derived. 

Our first observation of the above variational inequality based estimators is that the dual norm of Lagrange multiplier is relatively small, compared with the full error estimator. Indeed, we numerically tested above estimators for the three test cases, and present the representative result of Test case 3  in Table \ref{tab:comparison_nonsmooth_test3}, with $h=2^{-4}$. The dual norm is approximated in $V_{h/2}^{p}$, with the energy norm approximated by the dG energy norm.

\begin{table}[!htbp]
    \centering
    \begin{tabular}{|c|c|c|}
        \hline
        \multicolumn{2}{|c|}{$h=2^{-4}$}
        &  Ratio of Lagrange multiplier dual norm 
                      to total estimator \\
        \hline
        \multirow{2}{*}{$\mathbb{P}_1$}
        & $L^2$ corr. & $8.738\times10^{-4}$ \\
        \cline{2-3}
        & dG corr.   & $3.686\times10^{-4}$ \\
        \hline
        \multirow{2}{*}{$\mathbb{P}_2$}
        & $L^2$ corr. & $3.129\times10^{-3}$ \\
        \cline{2-3}
        & dG corr.   & $9.305\times10^{-4}$ \\
        \hline
    \end{tabular}
    \caption{Comparison of Lagrange multiplier based estimators on the nonsmooth Test case 3.}
    \label{tab:comparison_nonsmooth_test3}
\end{table}

\begin{figure}[!htbp]
    \centering
    \includegraphics[width=0.3\linewidth]{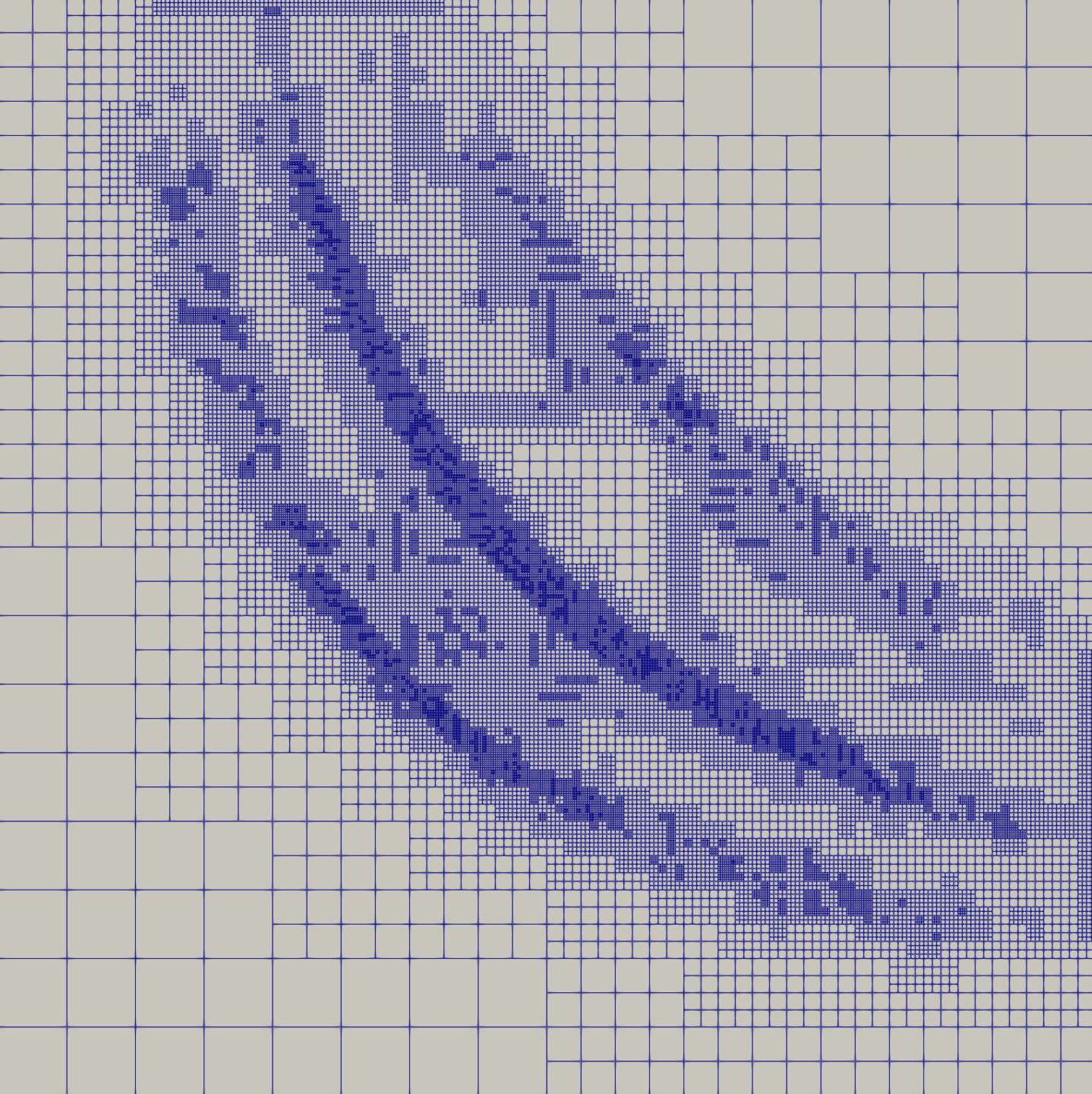}
    \includegraphics[width=0.3\linewidth]{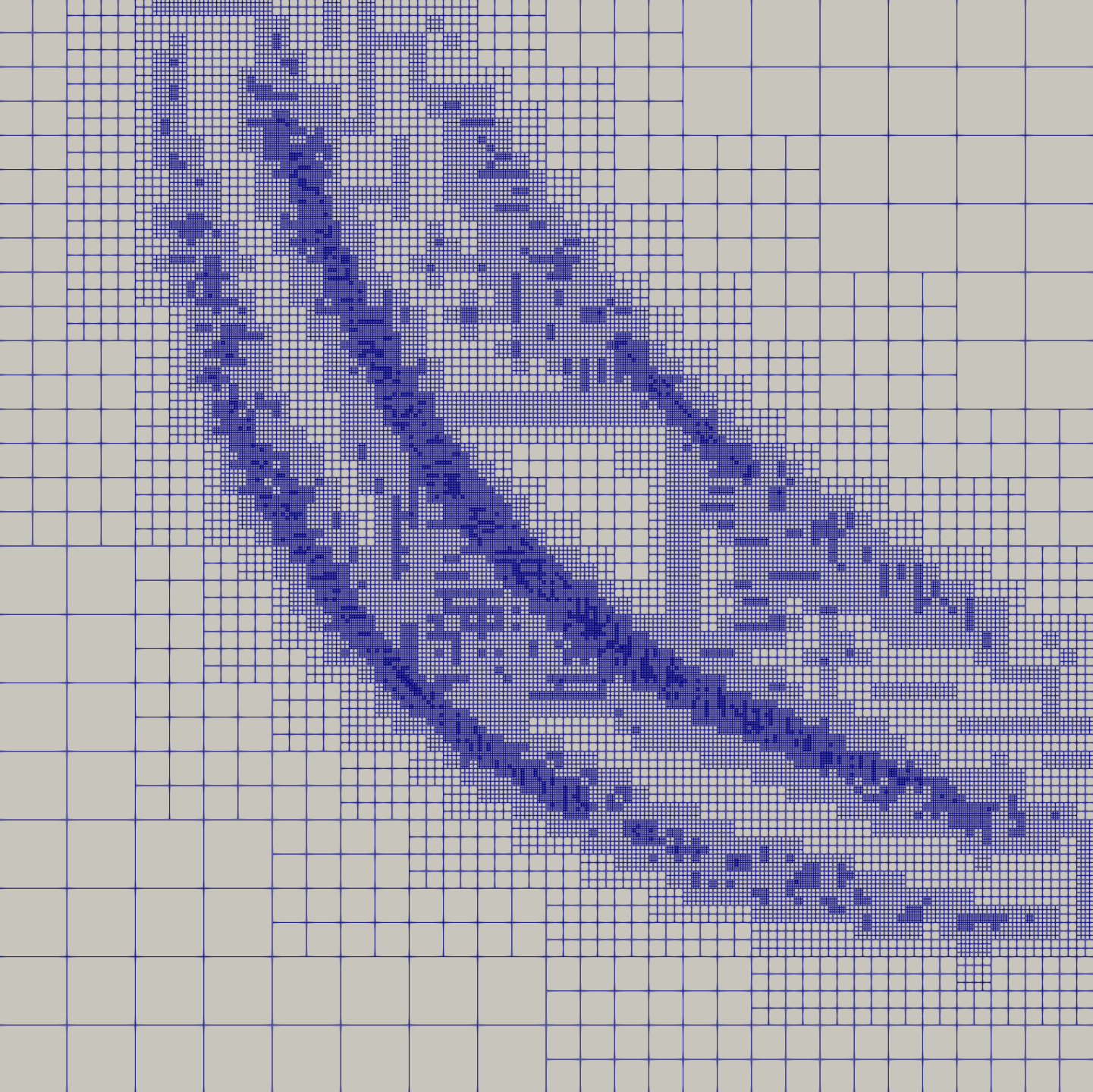}
    \includegraphics[width=0.3\linewidth]{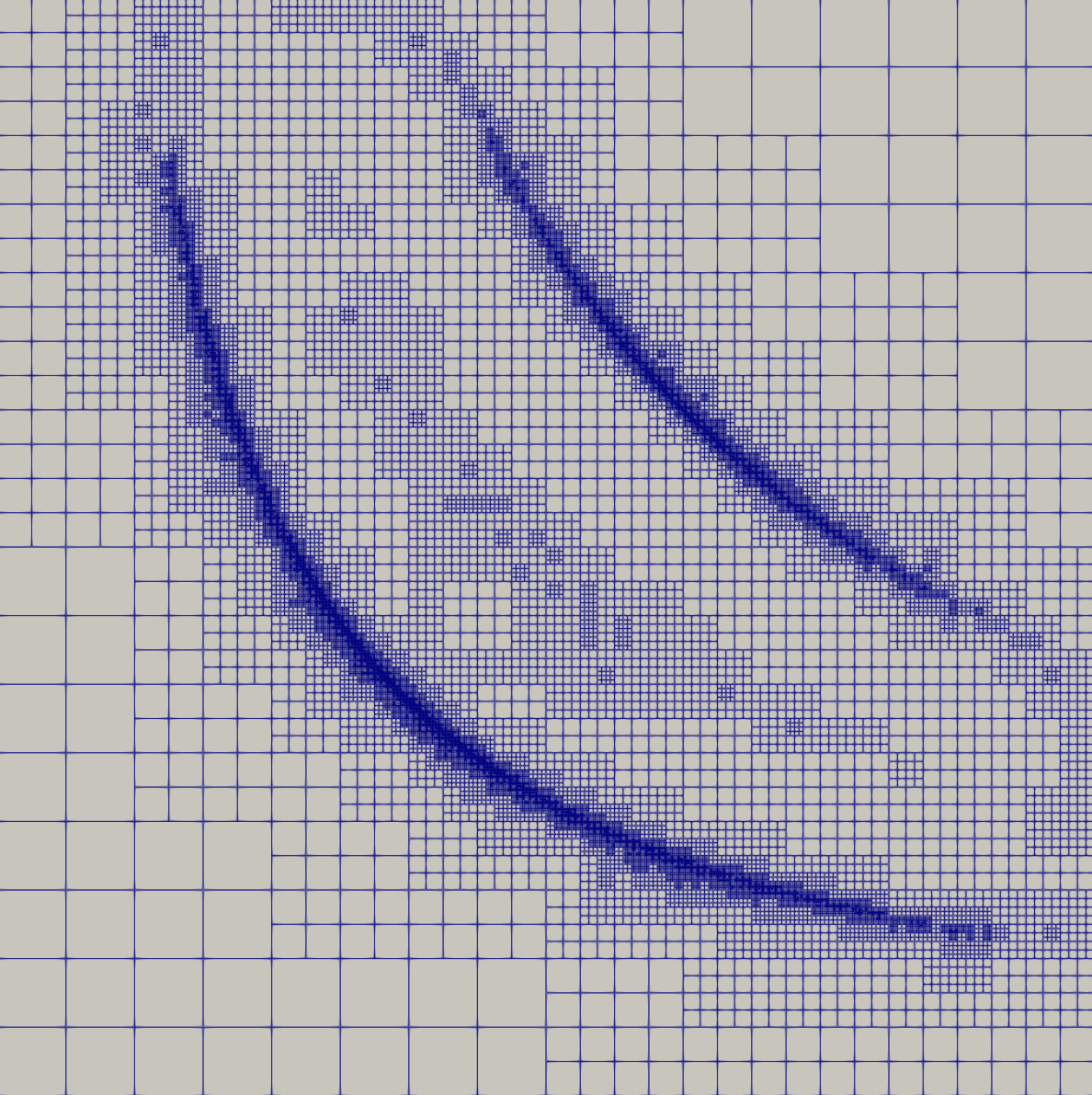}
    
    \caption{Meshes generated by adaptive scheme \eqref{eq_adaptive_scheme} with different corrections: left: $L^2$ correction \eqref{eq_scheme_step2_L2}; middle dG energy correction with Lagrange multiplier; right: $L^2$ correction with Lagrange multiplier.}
    \label{fig:adaptive_compare_LM}
\end{figure}

Then, we use the estimators \eqref{eq_estimator_LM_dG_energy} and \eqref{eq_estimator_LM_L2} (without dual norm) in the adaptive scheme \eqref{eq_adaptive_scheme}, and compare their behavior with our estimators. In Figure \ref{fig:adaptive_compare_LM}, we illustrate the final mesh with around $40000$ DOFs, produced by the adaptive scheme \eqref{eq_adaptive_scheme} with estimators \eqref{eq_eta_L2}, \eqref{eq_estimator_LM_dG_energy} (without dual norm) and \eqref{eq_estimator_LM_L2} (without dual norm), respectively. The parameters are the same as in Figure \ref{fig:adaptive_rate}. We note that only the $L^2$ correction based estimator is considered for our scheme here, since the dG energy correction based and prediction based estimators lead to similar results, as indicated in Figure \ref{fig:adaptive_ratio}.

From Figure \ref{fig:adaptive_compare_LM}, we observe that the Lagrange multiplier based estimator with dG correction produces a mesh close to the one delivered by our estimator which correctly refines cells with large error (compare with error plot in Figure \ref{fig:error_Kirby}), whereas the Lagrange multiplier based estimator with $L^2$ correction overestimates the error in bound-violating cells and consequently over-refines them. This suggests using our estimator instead of the ones used in the variational inequality literature.

\section{Concluding remarks}\label{sec_conclusion}
We developed and analyzed two prediction-correction dG methods for reaction-diffusion problems with anisotropic and heterogeneous diffusion. Both methods preserve prescribed nodal bounds and total mass. The energy correction retains the optimal dG-energy approximation order, while the $L^2$ correction retains the optimal $L^2$ order under the stated dual regularity assumption. Quasi-uniformity and the discrete inverse inequality are used for the additional optimal dG-energy estimate for the $L^2$ correction.

A key ingredient of the a priori analysis is the construction of a pointwise bound-preserving and cellwise mass-preserving comparison function. We also established an a posteriori transfer result that applies to any discrete correction of the predictor. It yields global reliability and efficiency up to data oscillation, while the local estimates contain an additional correction term.

The numerical experiments indicate that the energy correction produces smoother corrections in underresolved cases, whereas the $L^2$ correction is attractive when the relevant spatial features are resolved. 

In future work, we will consider applications of the current prediction-correction framework to nonlinear elliptic PDE, time-dependent PDE, as well as pointwise error estimates in the sense of \cite{Johnson_Schatz_Wahlbin_1987}.

\bibliographystyle{abbrvnat} 
\bibliography{bib}

\end{document}